\documentclass[11pt,reqno]{amsart}
\pdfoutput=1

\usepackage{amssymb,amscd,amsfonts,amsbsy,enumerate}
\usepackage{latexsym}
\usepackage{exscale}
\usepackage{amsmath,amsthm,amsfonts}
\usepackage{mathrsfs}
\usepackage{graphics}
\usepackage{esint}
\usepackage{enumitem} 
\usepackage{tikz}
\usepackage{pgfplots}
\pgfplotsset{compat=1.10}
\usepgfplotslibrary{fillbetween}
\usetikzlibrary{patterns}
\usepackage[colorlinks=true, pdfstartview=FitV, linkcolor=blue, citecolor=red, urlcolor=blue, pagebackref]{hyperref}

\usepackage{color}

\numberwithin{equation}{section}

\allowdisplaybreaks

\newtheorem{theorem}{Theorem}[section]
\newtheorem{lemma}[theorem]{Lemma}
\newtheorem{corollary}[theorem]{Corollary}
\newtheorem{proposition}[theorem]{Proposition}
\newtheorem{definition}[theorem]{Definition}

\theoremstyle{remark}
\newtheorem{remark}[theorem]{Remark}

\newcommand{\rnn}{\mathbb{R}^{n+1}}
\newcommand{\rn}{\mathbb{R}^{n}}
\newcommand{\rnp}{\mathbb{R}^{n+1}_{+}}

\newcommand{\X}{\mathbf{X}}

\newcommand{\Y}{\mathbf{Y}}

\begin{document}
\allowdisplaybreaks

\title[Nonlinear Neumann problem]
{On a nonlinear Laplace equation related to the boundary Yamabe\\ problem in the upper-half space}

\author{Gael Diebou Yomgne}
\address{
Rheinische Friedrich-Wilhelms-Universit\"{a}t Bonn
\\
Endeniche Allee 60
\\
53115, Bonn, Germany} \email{gaeldieb@math.uni-bonn.de}

\thanks{The author acknowledges financial support 
from the DAAD through the program ''Graduate School Scholarship Programme, 2018'' (Number 57395813) and  the grant from the University of
Bonn (HCM)}

\date{\today}

\subjclass[2010]{35B07, 35B09, 35B65, 35B33, 35C06, 35C15,42B37.}

\keywords{Boundary Yamabe problem, scaling invariance, Neumann data, Lorentz spaces, positive solutions, critical exponent}

\begin{abstract}
We consider in this paper the nonlinear elliptic equation with Neumann boundary condition 
\begin{align*}
\begin{cases}
\Delta u=a|u|^{m-1}u\,\,\mbox{ in }\,\,\rnp\\
\dfrac{\partial u}{\partial t}=b|u|^{\eta-1}u+f\,\,\mbox{ on }\,\,\partial\rnp.
\end{cases}
\end{align*}  
For $a,b\neq 0$, $m>\frac{n+1}{n-1}$, $(n>1)$, $\eta=\frac{m+1}{2}$ and small data $f\in L^{\frac{nq}{n+1},\infty}(\partial\rnp)$, $q=\frac{(n+1)(m-1)}{m+1}$ we prove that the problem is solvable. More precisely, we establish existence, uniqueness and continuous dependence of solutions on the boundary data $f$ in the function space $\X^{q}_{\infty}$ where
\[\|u\|_{\X^{q}_{\infty}}=\sup_{t>0}t^{\frac{n+1}{q}-1}\|u(\cdot,t)\|_{L^{\infty}(\partial\rnp)}+\|u\|_{L^{\frac{q(m+1)}{2},\infty}(\rnp)}+\|\nabla u\|_{L^{q,\infty}(\rnp)}.
\] 
As a direct consequence, we obtain the local regularity property $C^{1,\nu}_{loc}$, $\nu\in (0,1)$ of these solutions as well as energy estimates for certain values of $m$. Boundary values decaying faster than  $|x|^{-(m+1)/(m-1)}$, $x\in \rn\setminus\{0\}$ yield solvability and this decay property is shown to be sharp for positive nonlinearities.

Moreover, we are able to show that solutions inherit qualitative features of the boundary data such as positivity, rotational symmetry with respect to the $(n+1)$-axis, radial monotonicity in the tangential variable and homogeneity. When $a,b>0$, the critical exponent $m_c$ for the existence of positive solutions is identified, $m_c=(n+1)/(n-1)$. 
\end{abstract}

\maketitle

\allowdisplaybreaks

\normalsize

\section{Introduction and results}
\setcounter{equation}{0}
\label{S-1}

\label{S-1}
The classical Yamabe problem \cite{Tru} consists of finding a metric with constant scalar curvature which is conformal to a Riemannian metric of a closed manifold. There is a natural extension of this problem to compact manifolds with boundary due to Escobar \cite{Esco-1,Esco-2}. Let $(M,g)$ be a compact Riemannian manifold with dimension $n>2$ and $\partial M\neq \emptyset$, can one find a metric $\overline{g}$ in the conformal class of $g$ such that $R_{\overline{g}}$ and $\kappa_{\overline{g}}$ are both constant where $R_{\overline{g}}$ and $\kappa_{\overline{g}}$ are the scalar curvature of $\overline{g}$ and the mean curvature of the boundary $\partial M$, respectively?\\ Let  
\[\rnp:=\big\{(x,t), \hspace{0.1cm}x=(x_1,x_2,...,x_n)\in \rn,\hspace{0.1cm}t>0\big\}\] 
be the higher dimensional upper-half space and $B_1(0)$ the unit ball in $\rnn$. Since $M=B_1(0)$ is conformally equivalent to $\rnp$, the corresponding boundary Yamabe problem is equivalent to looking for positive classical solutions of the boundary value problem
\begin{align}\label{eq:Yam-eq}
\begin{cases}
\Delta u=u^{2^{\ast}-1}\,\,\mbox{ in }\,\,\rnp\\
\dfrac{\partial u}{\partial t}=u^{\eta^{\ast}-1}\,\,\mbox{ on }\,\,\partial\rnp
\end{cases}
\end{align}
with $2^{\ast}=\frac{2(n+1)}{n-1}$ and $\eta^{\ast}=\frac{2n}{n-1}$ which appear to be  the critical exponents for the Sobolev embedding $\mathscr{D}^{1,2}(\rnp)\hookrightarrow L^{2^{\ast}}(\rnp)$ and the trace principle $\mathscr{D}^{1,2}(\rnp)\hookrightarrow L^{\eta^{\ast}}(\partial\rnp)$,  respectively. Here, $\mathscr{D}^{1,2}(\rnp)$ denotes the space of restrictions to $\rnp$ of elements in $\mathscr{D}^{1,2}(\rnn)$, the completion of the space of smooth compactly supported functions in $\rnn$ with respect to the homogeneous Sobolev norm $\|\nabla u\|_{L^{2}(\rnn)}$.\\ In this paper, we are interested in the following general model
\begin{align}\label{eq:main-eq}
\begin{cases}
\Delta u=a|u|^{m-1}u\,\,\mbox{ in }\,\,\rnp\\
\dfrac{\partial u}{\partial t}=b|u|^{\eta-1}u+f\,\,\mbox{ on }\,\,\partial\rnp
\end{cases}
\end{align}
where $n>1$, $m,\eta>1$ and $a,b\neq 0$ are two real constants. Given a function $f$ defined on the boundary and assuming a minimal smoothness requirement, our goal in this paper is twofold. Firstly, we would like to develop a solvability theory for Eq. \eqref{eq:main-eq} and study the qualitative features of the corresponding solutions. Secondly, we aim at finding optimal conditions on the exponents and boundary data for which this well-posedness remains valid.\\

We briefly review a few known results on problem \eqref{eq:main-eq}. For $a\neq 0$, $b=0$ and zero boundary data $f=0$, Eq. \eqref{eq:main-eq} corresponds to the Neumann problem for the standard Lane-Emden equation. Its solvability highly depends on the exponent $m$ and the underlying geometry. Indeed, when posed in a  domain $\Omega\subset\rnn$, $(n>1)$ with smooth boundary, a weak solution may be sought for using a variational approach (\cite{S,QS}) under the condition $m\leq M_c$ so that the continuous embedding $\mathscr{D}^{1,2}(\Omega)\hookrightarrow  L^{m+1}(\Omega)$ holds and is compact provided $m<M_c=(n+3)/(n-1)$ and $\Omega$ bounded. In the critical case $m=M_c$, infinitely many solutions exist if $\Omega$ is a ball, see \cite{CK} and references therein for other interesting results. Also, it is known that any weak solution in bounded domain is $C^{2}$ regular provided $m\leq M_c$; this condition being sharp as counterexamples show \cite{DDF}. The study of harmonic functions subject to nonlinear Neumann data ($a=0$, $b\neq 0$) has also attracted a lot of attention in recent years. When $f$ belongs to the Lebesgue space $L^{1}(\Omega)$, with $\Omega$ bounded and $C^{2}$, existence of very weak solutions in the subcritical regime  $\eta<n/(n-1)$ and their regularity was analyzed in \cite{QR} using layer potentials techniques introduced in \cite{Fabes}. The exponent $\eta_{c}=n/(n-1)$ was shown to be critical for existence of bounded solutions. We point out, however, that critical exponent is a notion  depending on the regularity of the domain. The supercritical case $\eta>\eta_c$ in half-space was later investigated by other authors in \cite{FER-MEDEROS-MON}. Among  other results they obtained well-posedness whenever $f$ lies in a suitable Lebesgue space and has a small norm, by means of a nonvariational method.  Inspired by the above quoted works, a similar supercritical problem with singular potentials and data at the boundary was considered in \cite{AL1}. More general results of this flavour in connection to nonlocal equations is the main subject of the forthcoming work \cite{Yom}. In case $a$ and $b$ are simultaneously nonzero, the presence of the boundary datum in the problem causes the critical exponent to differ from that corresponding to the aforementioned scenarios, i.e. $a\neq 0$, $b=0$ or $b\neq 0$ and $f\equiv0$. A preliminary analysis using scaling may suggest the right candidate. In fact, $\eta=\frac{m+1}{2}$ is the unique exponent  so that the rescaled function 
\begin{equation}\label{scaling}
u_{\lambda}(x,t)=\lambda^{\frac{2}{m-1}}u(\lambda x,\lambda t)
\end{equation} solves  Eq. \eqref{eq:main-eq} with $f_{\lambda}(x)=\lambda^{\frac{m+1}{m-1}}f(\lambda x)$ in place of $f$ whenever $u$ is a classical solution. This suggests, for instance the choice of the datum $f$ from the Lorentz scale $L^{m^{\ast},\infty}(\rnp)$ where $m^{\ast}=\frac{n(m-1)}{m+1}$  since it captures the homogeneity of $f$. The condition $m^{\ast}>1$ which makes the latter space Banach is equivalent to the restriction $m>(n+1)/(n-1)$. We then anticipate solvability in this range, in particular the existence of positive solutions. Observe in passing that whenever this prediction is correct, the critical power for the boundary nonlinearity coincides with $\eta_c$.  
Regarding the study of positive solutions of Eq. \eqref{eq:main-eq}, it is clear that the signs of the constants $a$ and $b$ play a relevant role. As a matter of fact, when $f$ is identically zero on $\partial\rnp$, $m,\eta>1$ and $a,b>0$, it is long-established \cite{Lou-Zhu} that there is no positive classical solution. The remaining cases divided with respect to the signs of the constants is now fully understood, see the recent article \cite{Tang} and cited references therein. Now, we make precise the notion of solutions we study here. 
\begin{definition}\label{defn:int-eq}
Call $u$ a solution of Eq. \eqref{eq:main-eq} if $u$ satisfies the integral equation 
\begin{equation}\label{eq:int-eq}
u(x,t)=\mathscr{N}(b|u|^{\eta-1}u)(x,t)+\mathscr{G}(a|u|^{m-1}u)(x,t)+\mathscr{N}f(x,t),\quad(x,t)\in \rnp
\end{equation} 
where $\mathscr{G}$ is the Green potential for $\rnp$ associated to the Laplace operator (see Section \ref{S2}) and $\mathscr{N}$ is the convolution operator in $\rn$ with the Neumann kernel $\mathcal{K}_{t}$, $t>0$ defined by
\[\mathcal{K}_{t}(x):=\mathcal{K}(x,t)=\frac{\beta(n)}{(|x|^{2}+t^{2})^{\frac{n-1}{2}}};\hspace{0.2cm}x\in \rn,\hspace{0.1cm}\beta(n)=\pi^{\frac{n+1}{2}}\Gamma\bigg(\frac{n-1}{2}\bigg).\]
\end{definition}  

In order to make sense of the reformulation \eqref{eq:int-eq}, one would like to work in a functional framework whose norm is invariant with respect to the transformation \eqref{scaling}. This space should be defined in such a way that the nonlinear terms in the equation have correct mapping properties. Let  $1<q<n+1$ and set $q^{\ast}=\dfrac{(n+1)q}{n+1-q}$. We introduce the function space 
\begin{equation*}
\X^{q}_{\infty}=\mathscr{D}_{\infty}^{1,q}(\rnp)\cap L^{q^{\ast},\infty}(\rnp)\cap L_q^{\infty}(\rnp)   
\end{equation*} and equip it with the norm
\begin{equation}\label{norm}\|u\|_{\X^{q}_{\infty}}=\sup_{t>0}t^{\frac{n+1}{q}-1}\|u(\cdot,t)\|_{L^{\infty}(\partial\rnp)}+\|u\|_{L^{q^{\ast},\infty}(\rnp)}+\|\nabla u\|_{L^{q,\infty}(\rnp)}.
\end{equation}
Here, $\mathscr{D}_{\infty}^{1,q}(\rnp)=\{u\big|_{\rnp}: u\in \mathscr{D}_{\infty}^{1,q}(\rnn)\}$ where $\mathscr{D}_{\infty}^{1,q}(\rnn)$ is the completion of $C^{\infty}_0(\rnn)$ with respect to $\|\nabla u\|_{L^{q,\infty}(\rnn)}$ and 
$L^{\infty}_q(\rnp)$ collects all  measurable functions $u:\rnp\rightarrow  \mathbb{R}$ such that $\displaystyle\sup_{t>0}t^{\frac{n+1}{q}-1}\|u(\cdot,t)\|_{L^{\infty}(\rn)}$ is finite.
 With obvious modifications, we also consider the Lebesgue-based version of this function space which we denote by $\X^{q}$. Clearly, each of these framework carries out the structure of a complete normed space when endowed with their respective norms. Note that the latter gives the precise growth rate of solutions (and possibly that of its derivatives) near the boundary uniformly in the tangential variable. Moreover, the use of $\X^{q}_{\infty}$ allows one to obtain local regularity results directly from the solvability theory and also to prescribe boundary data with singularities.\\

 We prove (see Theorem \ref{thm:1}) that small data in  $L^{n(m+1)/(m-1),\infty}(\rn)$ and in particular any datum satisfying $|f(x)|\leq h(x)=\varepsilon_0|x|^{-(m+1)/(m-1)}$, $\varepsilon_0\lll 1$ gives rise to solvability via a fixed point argument. If $f$ rather decays slower than $h$ for large values of $|x|$, then we obtain nonexistence (Theorem \ref{thm:nonexistence-sol}). The proof relies on the implication: 
 \[u\hspace{0.1cm}\mbox{solves} \hspace{0.1cm} \eqref{eq:int-eq} \hspace{0.1cm}\mbox{in}\hspace{0.1cm} \X^{q}\hspace{0.1cm} \Longrightarrow\hspace{0.1cm} u\hspace{0.1cm}\mbox{ is a solution of}\hspace{0.1cm} Eq.\hspace{0.1cm} \eqref{eq:main-eq}\hspace{0.1cm} \mbox{in the sense of distributions}.\] We also construct solutions which are positive (for $a,b>0$), homogeneous, radially monotone in the tangential variable, rotationally symmetric around the vertical axis -- all these features once again are inherent to $f$. This is summarized in Theorem \ref{thm:4}. Finally, under some mild local assumptions on $f$, it is shown that Eq. \eqref{eq:int-eq} possesses no positive solutions when $m\in (1,(n+1)/(n-1)]$ (see Proposition \ref{prop:nonexistence-pos}). This reveals that $m_c=(n+1)/(n-1)$ is the sharp critical exponent. 

\subsection{Main results} The main theorems of this paper read as follows.
\begin{theorem}\label{thm:1} Let $a$ and $b$ be real constants not all simultaneously zero. Assume $m>\dfrac{n+1}{n-1}$, $n>1$ and put  $\eta=\dfrac{m+1}{2}$, $q=\dfrac{(n+1)(m-1)}{m+1}$. Then problem \eqref{eq:main-eq} is well-posed.\\ 
\textbf{$($i$)$} There exist $\varepsilon>0$ and $\vartheta:=\vartheta(\varepsilon)>0$ with the property that for every $f\in L^{\frac{nq}{n+1},\infty}(\rn)$ with $\|f\|_{L^{\frac{nq}{n+1},\infty}(\rn)}< \varepsilon$, there exists a function $u$ solution of Eq. \eqref{eq:int-eq} in $\X^{q}_{\infty}$ which is the only one satisfying the condition $\|u\|_{\X^{q}_{\infty}}\leq C\vartheta$ for some constant $C>0$.\\ 
\textbf{$($ii$)$} Moreover, let $m=1+\dfrac{4}{n}$ $($$n>2$$)$ and assume $f\in L^{\frac{2n}{n+2}}(\rn)$ with small norm. We have the energy-type estimate
\begin{equation}\label{energy-ineq}
 \bigg(\int_{\rnp}t|\nabla u|^{2}dX\bigg)^{1/2}\leq C\big(\|f\|_{L^{\frac{2n}{n+2}}(\rn)}+\|u\|^{\eta}_{\X^{q}}+\|u\|^{m}_{\X^{q}}\big)
\end{equation}
for some $C:=C(a,b)>0$, $a,b\neq 0$. Under the conditions $m=\dfrac{n+3}{n-1}$, there holds the energy estimate
\begin{equation}\label{energy-est-1}
\|\nabla u\|_{L^{2}(\rnp)}\leq C\big(\|f\|_{L^{2n/(n+1)}(\rn)}+\|u\|^{\eta}_{\X^{q}}+\|u\|^{m}_{\X^{q}}\big)
\end{equation} 
for some constant $C>0$ independent of $f$.
\end{theorem}

\vspace{0.51cm}
Our next result shows that the choice of $f$ in a smaller space together with the same smallness condition on $f$ in $L^{\frac{qn}{n+1},\infty}(\rn)$ produces a unique solution satisfying a higher integrability property. In more details, we have the following statement.
\begin{theorem}\label{thm:2}Let $a,b,m$ and $q$ falling under the scope of Theorem \ref{thm:1}. Let $1<p_0<n$ and suppose that $f\in L^{p_0,\infty}(\rn)\cap L^{\frac{nq}{n+1},\infty}(\rn)$. Then there exist $0<\varepsilon_{p_0}<\varepsilon$ and $\vartheta_{p_0}>0$ such that if $\|f\|_{L^{\frac{nq}{n+1},\infty}(\rn)}< \varepsilon_{p_0}$, then Eq. \eqref{eq:int-eq} is solvable in $\X^{\frac{(n+1)p_0}{n}}_{\infty}$ in the spirit of Theorem \ref{thm:1}.
\end{theorem}

\begin{corollary}[Regularity]\label{thm:3}
Granted the hypotheses of Theorem \ref{thm:1}, if $u$ is the  solution found, then   $u\in C^{1,\nu}_{loc}(\rn\times [t_0,\infty))$, $\nu\in (0,1)$ for any $t_0>0$. 
\end{corollary}   
This result is actually a consequence of Theorem \ref{thm:1}. For instance, it can be deduced from the observation that $u$ is bounded away from the origin uniformly in the tangential direction ($u$ has a membership to $\X^{q}_{\infty}$) and by a direct application of standard elliptic regularity theory. Indeed, since $u$ is bounded on $\rn\times [t_0,\infty)$ for any $t_0>0$, it follows that $u\in W^{2,p}(D)$, $1<p<\infty$ for any subset $D$ of $\rnp$ with $\overline{D}\subset \rn\times [t_0,\infty)$. Hence, one may choose $p$ sufficiently large ($p>n+1$) and use Sobolev inequality to obtain the desired conclusion. Actually, this regularity result can be strengthened by imposing a higher Sobolev  regularity on $f$ so as to obtain classical solutions.  
\begin{remark}\label{rmk:thm}It is easy to see that well-posedness for small data in $L^{\frac{n(m-1)}{m+1}}(\rn)$ immediately follows from (\textbf{i}). The following other observations can be made.
\begin{enumerate}
\item Our existence result shows in particular that solutions and their first order derivatives are such that \[|\nabla^{\kappa} u(x,t)|=O(t^{-\frac{2}{m-1}-\kappa}), \hspace{0.21cm} t\rightarrow \infty\hspace{0.1cm} \mbox{uniformly in}\hspace{0.1cm} x\in \rn; \hspace{0.1cm}\kappa=0,1.\] 
\item In case the underlying geometrical setting is $\Omega\subset \rn$, an open bounded domain of regularity $C^{1,\alpha}$, our well-posedness result remains valid for Lebesgue data in the natural candidate space where the norm of $L^{\infty}_q$ is now replaced by $\displaystyle\sup_{X\in \Omega}d(X)^{2/(m-1)}|u(X)|$ with  $d(X):=\mbox{dist} (X,\partial\Omega)$, the distance from $X\in \Omega$ to the boundary $\partial\Omega$.		
\item It is worth mentioning that our method allows one to also consider the  nonlinearity $N(u)=|X|^{-\gamma}u^{m}$, known as Hardy-type if $\gamma>0$ in which case it is necessary to require  $0<\gamma<\min(2,n+1)$ or H\'{e}non type whenever $\gamma<0$. 
\end{enumerate}
\end{remark}
We now turn our attention to the study of solutions enjoying more qualitative properties.
\begin{theorem}\label{thm:4}
Given the hypotheses of Theorem \ref{thm:1}, let $u$ be the solution constructed in $\X^{q}_{\infty}$. One has
\begin{enumerate}[label={($\textbf{a}_{\arabic*})$}]
\item (Positivity)\label{pos-sol} Assume $f$ is nonnegative in $\rn$ and positive in an open subset $\Omega\subset\rn$ with finite Lebesgue measure. If $a,b>0$ then $u>0$ in $\rnp.$  
\item(Rotational symmetry)\label{rot-sym}	Let $a,b\neq 0$. The function $f$ is radial in $\rn$ if and only if $u$ is symmetric with respect to any group of rotations in $\rnp$ leaving the vertical axis fixed.
\item(Radial monotonicity)\label{rad-mono} The solution $u$ (with $a,b\geq 0$) is radially nonincreasing in $x\in \rn$ if $f$ is so.
\item (Homogeneity)\label{self-similar} If $f$ is homogeneous of degree $-\frac{m+1}{m-1}$ and has a small norm in the Lorentz space $L^{\frac{n(m-1)}{m+1},\infty}(\rn)$, then the corresponding solution $u$ of \eqref{eq:int-eq} is homogeneous, that is, $u$ obeys \eqref{scaling}.  
\end{enumerate}
\end{theorem}

We also note that the remaining cases (e.g. $a$ and $b$ have opposite signs) are not covered by Theorem \ref{thm:4} whose proof essentially relies on successive iterations. 
Theorem \ref{thm:4}-\ref{pos-sol} tells us that the existence of signed solutions is a property inherited from the prescribed boundary value. However, it is well-known that Eq. \eqref{eq:main-eq} with $f\equiv0$ admits no positive classical solutions for $m,\eta>1$ and $a,b>0$, see \cite{Lou-Zhu}. On the other hand, when $a,b<0$, the conclusion remains the same for $0\leq m\leq \frac{n+3}{n-1}$, $-\infty<\eta\leq \frac{n+1}{n-1}$ as proved in \cite{YanYan-Zhang}.
A natural question to ask is: Do positive solutions exist if $m$ is confined in the range $\big(1,\frac{n+1}{n-1}\big]?$ The existence theory settled in Theorem \ref{eq:main-eq} already suggests $m_c=\frac{n+1}{n-1}$ as the exponent separating the regime of existence and nonexistence of positive solutions at least in the specific case $a,b>0$. Under a mild assumption on $f$, the next result shows that $m_c$ is indeed the critical exponent.
\begin{proposition}[Nonexistence]\label{prop:nonexistence-pos} Let $q$ as in Theorem \ref{thm:1} and $a,b>0$. Assume that $f$ is positive and bounded below near the origin. If $m\in \big(1,\frac{n+1}{n-1}\big)$, then Eq. \eqref{eq:int-eq} has no positive solutions in $\X^{q}_{\infty}$. Moreover, if $m=(n+1)/(n-1)$ and $f$ is bounded below locally, then no positive solutions exist in $\X^{q}_{\infty}$.  
\end{proposition}
Note that $\X^q$ is strictly contained in $\X^q_{\infty}$ so that the conclusions of the above proposition remains true in the former space. Essentially, the idea  consists of showing that the norm of the solution in $\X^q_{\infty}$ grows when $m$ is subcritical ($m<m_c)$. The proof also works in the case $a=0$ and the result reproduces a particular case of \cite[Theorem 1.10]{Yom} and equally solves the open question in \cite{AL1}.

\begin{remark}
The boundary datum $f$ in  Eq. \eqref{eq:main-eq} has the effect of shrinking the initial range (that of Eq. \eqref{eq:main-eq} with $f=0$) of nonexistence of positive solutions as depicted below (see Figure \ref{fig1}). 

 \begin{figure}[h]
 \centering	
\begin{tikzpicture} [scale=0.7]
\fill [blue,opacity=0.2, domain=1.4:17, variable=\x]
(1, 0)
-- (1,11)
-- plot[smooth,tension=.3] ({\x}, {(\x+3)/(\x-1)})
-- (17, 0)
-- cycle;
\fill [red,opacity=0.3, domain=1.2:17, variable=\x]
(1, 0)
-- (1,11)
-- plot[smooth,tension=.3] ({\x}, {(\x+1)/(\x-1)})
-- (17, 0)
-- cycle;
\draw [red,thin, domain=1.2:17, variable=\x] plot[smooth,tension=.3] ({\x}, {(\x+1)/(\x-1)});
\draw [blue,thin, domain=1.4:17, variable=\x] plot[smooth,tension=.3] ({\x}, {(\x+3)/(\x-1)});
\draw[thick,->] (0,0) -- (17,0);
\draw[thick,->] (0,0) -- (0,11);
\draw[thick,dashed] (1,0) -- (1,11);
\node[left,below] at (0,0) {$0$};
\node[left,below] at (1,0) {$1$};
\node[left,above] at (0,11) {$m$};
\node[right,below] at (17,0) {$n$};
\node[right,red] at (2,7.7) {$m_c=\frac{n+1}{n-1}$};
\node[right,blue] at (2,9) {$M_c=\frac{n+3}{n-1}$};
\end{tikzpicture}	
 \caption{The region below the critical curve $M_c(n)$ (resp. $m_c(n)$) indicates the nonexistence range relative to  Eq. \eqref{eq:main-eq} (with $f=0$ and $a,b\neq 0$ having same sign) (resp. for Eq. \eqref{eq:main-eq} with $a,b>0$ and $f\neq 0$).}\label{fig1}	
 \end{figure}
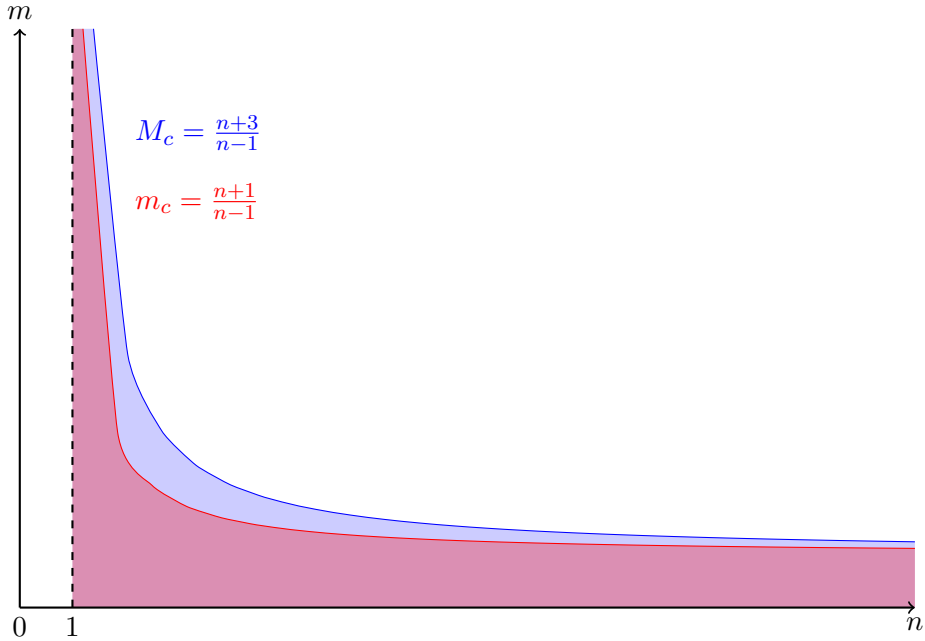	
\end{remark}
  
Moving on, we remark that Eq. \eqref{eq:int-eq} remains solvable for any datum $f$ defined on $\rn$, decaying faster than $c|x-x_0|^{-\frac{m+1}{m-1}}$, $x\neq x_0\in \rn$ for $|c|$ small. This decay condition seems to be optimal, at least for signed nonlinearities ($a,b>0$) as shown below. 
\begin{theorem}\label{thm:nonexistence-sol}
Take $f\in L^{1}_{loc}(\rn)$ and suppose that there exists $k\in [0,(m+1)/(m-1))$, $m>1$ such that $\displaystyle\liminf_{|x|\rightarrow \infty}(1+|x|^{k})|f(x)|\geq \varLambda$ for some constant $\varLambda>0$. Then Eq. \eqref{eq:int-eq} under the above considerations has no solution. 	
\end{theorem}     
A direct interpretation of this result is that boundary data which at infinity decay slower than the singularity $|x|^{-\frac{m+1}{m-1}}$, $x\in \rn\setminus\{0\}$ do not generate  solutions. The proof employs a test function method. This is possible to implement because  a solution to Eq. \eqref{eq:int-eq} is  a distributional solution of problem \eqref{eq:main-eq} under reasonable conditions. 
\section{Preliminaries}
\label{S2}
We start with standard notations and definitions of function spaces.  Subsequently, we state the auxiliary results which will find their usefulness later in Section \ref{S3}. \\
\subsection{Notations and definitions} 
Throughout the paper, $\mathbb{N}$ stands for the collection of all positive integers. We will identify $\rn$ with the boundary of the half-space $\rnp$, $\partial\rnp$ and $\rnn_{-}=\rnn\setminus \rnp$. A typical point in $\rnp$ will be denoted by $X=(x,t)$, $x\in \rn$ and $t>0$.\\
For $\alpha\in (0,n)$, the notation $I_{\alpha}$ is reserved for the Riesz potential of order $\alpha$, that is, the convolution operator with singular kernel $c_{\alpha,n}|x|^{\alpha-n}$, $c_{\alpha,n}=\pi^{-\frac{n}{2}}2^{-(\alpha+1)}\frac{\Gamma(n-\alpha)}{\Gamma(\alpha/2)}$ which transforms a Schwartz function $f$ in $\rn$ via Fourier transform according to \[\widehat{I_{\alpha}f}(\xi)=\widehat{(-\Delta)^{-\frac{\alpha}{2}}f}(\xi)=|\xi|^{-\alpha}\widehat{f}(\xi).\]
Given a set $A\subset \rn$, denote by $\textbf{1}_{A}$ its characteristic function and by $|A|$ its Lebesgue measure. Let  $f:\rn\rightarrow \mathbb{R}$ be a measurable function and $d_{f}(\tau)=|\{|f|> \tau\}|$ its distribution function. Call $f^{\ast}:[0,\infty)\rightarrow [0,\infty)$ the decreasing rearrangement of $f$ defined by \[f^{\ast}(\lambda)=\inf\{\tau>0:d_{f}(\tau)\leq \lambda\}\] and let $f^{\ast\ast}:(0,\infty)\rightarrow [0,\infty)$ be the associated maximal function given by $f^{\ast\ast}(\tau)=\frac{1}{\tau}\int_{0}^{\tau}f^{\ast}(s)ds$. For $1<p<\infty$, the Lorentz space  $L^{p,q}(\rn)$ is defined as 
	\[L^{p,q}(\rn)=\bigg\{f:\rn\rightarrow \mathbb{R} \hspace{0.1cm}\text{measurable} : \|f\|^{q}_{L^{p,q}(\rn)}=\int_{0}^{\infty}[s^{1/p}f^{\ast\ast}(s)]^{q}\frac{ds}{s}<\infty\bigg\},\hspace{0.12cm} q\geq 1
	\]
	\[\hspace{-2.2cm}L^{p,\infty}(\rn)=\bigg\{f:\rn\rightarrow \mathbb{R} \hspace{0.1cm}\text{measurable} : \|f\|_{L^{p,\infty}(\rn)}=\sup_{0<s<\infty}s^{1/p}f^{\ast\ast}(s)<\infty\bigg\}.
	\]
These spaces may alternatively be defined by means of $d_{f}$ using the quasi-norm 
\[|u|_{L^{p,q}(\rn)}=\bigg(p\int_{0}^{\infty}[sd_{f}(s)^{1/p}]^{q}\dfrac{ds}{s}\bigg)^{1/q}\]
with usual modification when $q=\infty$. Note that the functionals $|u|_{L^{p,q}}$ and $\|u\|_{L^{p,q}}$ are equivalent and in case $q=\infty$, another useful equivalent norm exists namely, 
\begin{equation*}
\|u\|_{L^{p,\infty}(\rn)}\approx\sup_{\Omega\subset \rn}|\Omega|^{1/p-1}\int_{\Omega}|u(y)|dy    
\end{equation*}
where the supremum is taken over all open subsets of $\rn$.
 A standard real interpolation result in Lorentz spaces reads as follows, see e.g. \cite{BeLo}. Assume $1<p_1<p_2<\infty$, $1\leq q,q_1,q_2\leq \infty$, one has  \[
[L^{p_1,q_1}(\rn),L^{p_2,q_2}(\rn)]_{\vartheta,q}=L^{p,q}(\rn),\quad \frac{1}{p}=\frac{1-\vartheta}{p_1}+\frac{\vartheta}{p_2}
\]
where $\vartheta\in (0,1)$. Recall the generalized H\"{o}lder inequality in the framework of Lorentz spaces due to O'Neil \cite{O}.	Let $p_1,p_2\in (1,\infty)$ and $1\leq q_1,q,q_2\leq \infty$ such that $\frac{1}{q}\leq  \frac{1}{q_1}+\frac{1}{q_2}$. For all $f\in L^{p_1,q_1}(\rn)$ and $g\in L^{p_2,q_2}(\rn)$, we have $fg\in  L^{p,q}(\rn)$ and there holds 
\begin{equation}\label{eq:Holder-ineq}
\|fg\|_{ L^{p,q}(\rn)}\leq C\|f\|_{L^{p_1,q_1}(\rn)}\|g\|_{L^{p_2,q_2}(\rn)}
\end{equation}
	provided $\frac{1}{p}= \frac{1}{p_1}+\frac{1}{p_2}$. 
\subsection{Harmonic functions with Neumann data} Consider the boundary value problem
\begin{align}\label{eq:Neumann-eq}
\begin{cases}
\Delta v=0\,\,\mbox{ in }\,\,\rnp\\
\dfrac{\partial v}{\partial\eta}=f\,\,\mbox{ on }\,\,\partial\rnp
\end{cases}
\end{align}
where  $f$ is continuous and belongs to the weighted Lebesgue space $L^1\big(\rn,(1+|x|^{2})^{-\frac{n-1}{2}}dx\big)$. Then $v=\mathscr{N}f$ with $\mathscr{N}$ as in Definition \ref{defn:int-eq} solves Eq. \eqref{eq:Neumann-eq}, see for instance \cite{Arm}. In the sequel, we establish key estimates for the operator $\mathscr{N}$ in the setting of Lorentz spaces.
\begin{proposition}\label{prop:linear-bounds}
Let  $0< \mu\leq \infty$. Given $f$ in $L^{p,\mu}(\rn)$, there exists  $C_1:=C_1(n,p,\mu)>0$ such that
\begin{equation}\label{eq:est1}
\|\nabla^{\kappa}\mathscr{N}f(\cdot,t)\|_{L^{q,\mu}(\rn)}\leq C_1t^{-n(\frac{1}{p}-\frac{1}{q})+1-\kappa}\|f\|_{L^{p,\mu}(\rn)}\hspace{0.1cm}\forall\hspace{0.1cm} t>0,\hspace{0.1cm}\frac{np}{n+(\kappa-1)p}<q\leq \infty; \hspace{0.1cm}\kappa=0,1
\end{equation}
for $1<p<\infty$ with $p<n$ when $\kappa=0$. Moreover, we have
\begin{equation}\label{eq:est2}
\|\mathscr{N}g\|_{L^{r_1,\mu}(\rnp)}\leq C_2\|g\|_{L^{\ell_1,\mu}(\rn)}, \hspace{0.2cm}\ell_1=\frac{nr_1}{n+1+r_1},\hspace{0.2cm} \frac{n+1}{n-1}<r_1<\infty
\end{equation}
for all $g\in L^{\ell_1,\mu}(\rn)$ and for some constant  $C_2:=C_2(n,\ell_1,\mu)>0$  and 
\begin{equation}\label{eq:est2'}
\|\nabla\mathscr{N}g\|_{L^{r_2,\mu}(\rnp)}\leq C_3(n,\ell_2,\mu)\|g\|_{L^{\ell_2,\mu}(\rn)}, \hspace{0.2cm}\ell_2=\frac{nr_2}{n+1},\hspace{0.2cm} \frac{n+1}{n}<r_2<\infty.
\end{equation}
for every $g\in L^{\ell_2,\mu}(\rn)$.
\end{proposition}
\begin{proof}
We easily verify that $|\nabla^{\kappa}\mathcal{K}(x,t)|\leq C(|x|^{2}+t^2)^{-\frac{n-1+\kappa}{2}}$, $\kappa=0,1$ and for $t>0$ fixed, \[|\nabla^{\kappa}\mathcal{K}(\cdot,t)|\in L^{\varLambda}(\rn)\hspace{0.2cm} \mbox{provided}\hspace{0.2cm} (n-1+\kappa)\varLambda>n.\] Also, for $f\in L^p(\rn)$, $\nabla^{\kappa}v(x,t)=\nabla^{\kappa}\mathscr{N}f(x,t)=(\nabla^{\kappa}\mathcal{K}(\cdot,t)\ast f)(x)$ and should be understood in the sense of distributions when $\kappa=1$. By using Young's convolution inequality, one finds that 
\begin{align}\label{identity}
\nonumber\|\nabla^{\kappa}v(\cdot,t)\|_{L^{q}(\rn)}&\leq C\|\nabla^{\kappa}\mathcal{K}(\cdot,t)\|_{L^{\varLambda}(\rn)}\|f\|_{L^{p}(\rn)}\\
&\leq Ct^{-(n-1+\kappa)}t^{n/\varLambda}\big\|(1+|\cdot|^{2})^{-\frac{n-1+\kappa}{2}}
\big\|_{L^{\varLambda}(\rn)}\|f\|_{L^{p}(\rn)}
\end{align}
where $1/p+1/\varLambda>1$ (which implies $1/p>(1-\kappa)/n$) and $\dfrac{1}{q}=\dfrac{1}{p}+\dfrac{1}{\varLambda}-1$. This leads to the estimate 
\begin{align*}
\|\nabla^{\kappa}v(\cdot,t)\|_{L^{q}(\rn)}&\leq Ct^{-n(\frac{1}{p}-\frac{1}{q})+1-\kappa}\|f\|_{L^{p}(\rn)}
\end{align*}	
which after real interpolation gives the desired bound \eqref{eq:est1}. One may establish the third estimate \eqref{eq:est2} as follows. Let $r_1,\ell_1\in (1,\infty)$, take $g\in L^{\ell_1,\mu}(\rn)$ and observe that \[\|\mathscr{N}g\|_{L^{r_1}(\rnp)}= \big\|\|\mathscr{N}g(x,\cdot)\|_{L^{r_1}(\mathbb{R}_{+},dt)}\big\|_{L^{r_1}(\rn)}.\] Next, using Minkowski's inequality, we estimate the inner norm   
\begin{align*}
\|\mathscr{N}g(x,\cdot)\|_{L^{r_1}(0,\infty)}&\leq \int_{\mathbb{R}^{n}}\|\mathcal{K}(x-y,\cdot)\|_{L^{r_1}(0,\infty)}|g(y)|dy\\
&\leq C\int_{\mathbb{R}^{n}}\|\mathcal{K}(x-y,\cdot)\|_{L^{r_1}(0,\infty)}|g(y)|dy\\
&\leq C\int_{\mathbb{R}^{n}}\bigg(\int_{0}^{\infty}\dfrac{dt}{(|x-y|^2+t^2)^{\frac{(n-1)r_1}{2}}}\bigg)^{1/r_1}|g(y)|dy\\
&\leq C\int_{\mathbb{R}^{n}}\dfrac{|g(y)|dy}{|x-y|^{n-(1+1/r_1)}}\\
&\leq CI_{\frac{1}{r_1}+1}|g|(x).
\end{align*} 
As such, we deduce from \eqref{identity} and the mapping properties of the Riesz potential in Lebesgue spaces (see e.g. \cite{St}) that
\begin{align*}
\|\mathscr{N}g\|_{L^{r_1}(\rnp)}&\leq C\big\|I_{1+\frac{1}{r_1}}|g|\big\|_{L^{r_1}(\rn)}\leq C\|g\|_{L^{\ell_1}(\rn)}
\end{align*}
provided  $\dfrac{1}{r_1}=\dfrac{1}{\ell_1}-\dfrac{1+\frac{1}{r_1}}{n}$. By real interpolation  we get \eqref{eq:est2}. Finally, regarding the gradient estimate \eqref{eq:est2'},  one may follow the steps of the above proof making use of the pointwise gradient bound on the Neumann kernel $\mathcal{K}$. This achieves the proof of Proposition \ref{prop:linear-bounds}.
\end{proof}   

\subsection{Poisson equation and Nonlinear estimates}
In this part, we are interested in the mapping properties of the nonlinearity which may be easily deduced from the study of an appropriate inhomogeneous problem. Recall the space $\X^{q}_{\infty}$ introduced in Section \ref{S-1}. Here we consider an analogous space which we denote by $\Y^{p}_{\infty}$, $1<p<\infty$ with norm \[\|F\|_{\Y^{p}_{\infty}}:=\sup_{t>0}t^{\frac{n+1}{p}}\|F(\cdot,t)\|_{L^{\infty}(\partial\rnp)}+\|F\|_{L^{p,\infty}(\rnp)}.\]
Likewise, we write $\Y^p$ to denote the analogous space with norm defined using the Lebesgue norm instead. 
Define the Green potential $\mathscr{G}$ via the representation formula \[\mathscr{G}F(X)=\int_{\rnp}G(X,Y)F(Y)dY\] for a suitable function $F$ where $G(\cdot,\cdot)$ is the Green kernel for the Laplace operator in $\rnp$ with Neumann boundary condition, that is,
\begin{align}\label{eq:green-eq}
\begin{cases}
\Delta G(X,\cdot)=\delta_{X}(\cdot)\,\,\mbox{ in }\,\,\rnp\\
\dfrac{\partial G(X,\cdot)}{\partial t}=0\,\,\mbox{ on }\,\,\partial\rnp
\end{cases}
\end{align}
in the sense of distributions where $\delta_{X}$ is the Dirac's distribution with mass at $X$. Thus, $\mathscr{G}F$ formally solves the boundary value problem $\Delta \mathscr{G}F=F$ in $\rnp$, $\partial_t\mathscr{G}F\big|_{\partial\rnp}=0$ and  $G(\cdot,\cdot)$ is a positive function which assumes an explicit form given by (see e.g. \cite{CP}) 
\[G(X,Y)=\gamma_n\bigg[|X-Y|^{-(n-1)}+|X-Y^{\ast}|^{-(n-1)}\bigg],\quad \gamma_n=\frac{1}{(n-1)\sigma_{n+1}} \]
for $X\in \overline{\rnp}$, $Y\in \rnp$,  $X\neq Y$ where $\sigma_{n+1}$ is the surface area of the unit sphere of $\rnn$ and
$Y^{\ast}=(y_1,...,-y_{n+1})$ is the reflection of the point $Y$ across the hyperplane $\{y_{n+1}=0\}$. Since $|X-Y|\leq |X-Y^{\ast}|$, it follows that $G$ satisfies the pointwise estimate 
\begin{equation}\label{eq:pointwise-bound}
G(X,Y)\leq 2\gamma_n|X-Y|^{1-n}.
\end{equation}
Moreover, we also have the gradient bound 
\begin{equation}\label{eq:grad-bound}
|\nabla_XG(X,Y)|\leq C|X-Y|^{-n}.
\end{equation}
In light of \eqref{eq:pointwise-bound} and \eqref{eq:grad-bound} we record the following bounds on the Green potential.   
\begin{lemma}\label{lem:Green-pot-bound}Let $n>1$ and $1<p<\dfrac{n+1}{2}$. Set $q=\dfrac{(n+1)p}{n+1-p}$. The Green potential $\mathscr{G}:\Y^{p}_{\infty}\rightarrow\X^{q}_{\infty}$  $($resp. $\mathscr{G}:\Y^{p}\rightarrow\X^{q}$$)$ continuously with the accompanying estimates
\begin{equation}\label{eq:nonlin-bound}
\|\mathscr{G}F\|_{\X^{q}_{\infty}}\leq C\|F\|_{\Y^{p}_{\infty}}\hspace{0.2cm}\mbox{for all}\hspace{0.2cm}F\in \Y^{p}_{\infty},
\end{equation}
and 	
\begin{equation}\label{eq:nonlinear-bd-Lebesgue}
\|\mathscr{G}F\|_{\X^{q}}\leq C\|F\|_{\Y^{p}}\hspace{0.2cm}\mbox{for every}\hspace{0.2cm}F\in \Y^{p} 
\end{equation}
for some constant $C:=C(n,p)>0$ independent of $F$. Furthermore, we have \begin{equation}\label{eq:grad-est-Green-kernel}
\sup_{t>0}t^{\frac{n+1}{p}-1}\|\nabla \mathscr{G}F(\cdot,t)\|_{L^{\infty}(\partial\rnp)}\leq C\|F\|_{L^{p}(\rnp)}.
\end{equation}
\end{lemma}
Note that if $F\in L^{p,\infty}(\rnp)$, then a version of  \eqref{eq:grad-est-Green-kernel} holds with the $L^{p,\infty}$-norm on R.H.S. We only present the proof of \eqref{eq:nonlin-bound} and \eqref{eq:grad-est-Green-kernel}. That of the estimate \eqref{eq:nonlinear-bd-Lebesgue} is essentially obtained via similar arguments.
\begin{proof}[Proof of Lemma \ref{lem:Green-pot-bound}]
Let $2p<n+1$ and $q>1$ with $\dfrac{n+1}{q}=\dfrac{n+1}{p}-1$. Assume that $F\in \Y^{p}_{\infty}$. 
 We first show that
\begin{equation}\label{eq:weighted-sup}
\displaystyle\sup_{t>0}t^{\frac{n+1}{p}-2}\|\mathscr{G} F(\cdot,t)\|_{L^{\infty}(\rn)}\leq \|F\|_{\Y^{p}_{\infty}}.
\end{equation}
Fix $(x,t)\in \rnp$ and let $B_t(x)\subset \rn$ denote the closed ball of radius $t>0$ with center at $x\in \rn$. Decompose the integral 
\begin{equation*}
\mathscr{G}F(x,t)=\int_{\rnp}G(x,t,Y)F(Y)dY=A_1+A_2+A_3+A_4
\end{equation*}
where		
\[A_1=\int_{B_{t}(x)}\int_{0}^{t/2} G(x,t,Y)F(Y)dY,\hspace{0.1cm}A_2=\int_{B_{t}(x)}\int_{t/2}^{2t} G(x,t,Y)F(Y)dY,\]
\[A_3=\int_{\rn\setminus B_{t}(x)}\int_{0}^{2t} G(x,t,Y)F(Y)dY,\hspace{0.1cm}A_4=\int_{\rn}\int_{2t}^{\infty} G(x,t,Y)F(Y)dY.\]
Next, we estimate each of these integrals using the pointwise bound \eqref{eq:pointwise-bound} repeatedly. Starting with $A_1$, we have
\begin{align*}
|A_1|&\leq\int_{B_{t}(x)}\int_{0}^{t/2}G(x,t,Y)|F(Y)|dY\\
&\leq C\int_{B_{t}(x)}\int_{0}^{t/2}\frac{|F(Y)|}{(|x-y|^{2}+(t-y_{n+1})^{2})^{\frac{n-1}{2}}}dy_{n+1}dy\\
&\leq C\|F\|_{L^{p,\infty}(B_{t}(x)\times (0,\frac{t}{2}))}\bigg\|\big[|x-y|^{2}+(t-y_{n+1})^{2}\big]^{\frac{1-n}{2}}\bigg\|_{L^{p',1}(B_{t}(x)\times (0,\frac{t}{2}))}\\
&\leq C\|F\|_{L^{p,\infty}(\rnp)}t^{1-n}t^{\frac{n+1}{p'}}\\
&\leq Ct^{2-\frac{n+1}{p}}\|F\|_{\Y^{p}_{\infty}}
\end{align*}
where $1/p+1/p'=1$. Note that we have utilized the generalized H\"{o}lder's inequality \eqref{eq:Holder-ineq} in order to derive the third bound in the above chain of estimates. On the other hand, one has
\begin{align*}
|A_2|&\leq \int_{B_{t}(x)}\int_{t/2}^{2t}G(x,t,Y)|F(Y)|dY\\
&\leq Ct^{-\frac{n+1}{p}}\sup_{y_{n+1}>0}y_{n+1}^{\frac{n+1}{p}}\|F(\cdot,y_{n+1})\|_{L^{\infty}(\rn)}\int_{B_{t}(x)}\int_{t/2}^{2t} G(X,Y)dY\\
&\leq C\sup_{y_{n+1}>0}y_{n+1}^{\frac{n+1}{p}}\|F(\cdot,y_{n+1})\|_{L^{\infty}(\rn)}t^{-\frac{n+1}{p}}t\int_{B_{t}(x)} |x-y|^{1-n}dy\\
&\leq Ct^{2-\frac{n+1}{p}}\|F\|_{\Y^{p}_{\infty}}.
\end{align*}
By invoking the generalized H\"{o}lder's inequality once again,  we arrive at
\begin{align*}
|A_3|&= \bigg|\int_{\rn\setminus B_{t}(x)}\int_{0}^{2t}G(X,Y)F(Y)dY\bigg|\\
&\leq C\int_{\rn\setminus B_{t}(x)}\int_{0}^{2t}|X-Y|^{1-n}|F(Y)|dY\\
&\leq C\sum_{i=1}^{\infty}\int_{2^{i}B_{t}(x)\setminus 2^{i-1}B_{t}(x)}\int_{0}^{2t}|X-Y|^{1-n}|F(Y)|dY\\
&\leq C\sum_{i=1}^{\infty}\int_{2^{i}B_{t}(x)\setminus 2^{i-1}B_{t}(x)}\int_{0}^{2t}|x-y|^{1-n}|F(Y)|dY\\
&\leq C\|F\|_{L^{p,\infty}(\rnp)}t^{1-n}\sum_{i=1}^{\infty}2^{-i(n-1)}(2^{i}t)^{\frac{n+1}{p'}}\\
&\leq Ct^{2-\frac{n+1}{p}}\|F\|_{\Y^{p}_{\infty}}
\end{align*}
because $n+1>2p$. Finally, with the same ingredients, we bound $A_4$ as follows: 
\begin{align}\label{A4}
\nonumber|A_4|&\leq \int_{\rn}\int_{2t}^{\infty}G(X,Y)|F(Y)|dY\\
\nonumber&\leq C\int_{\rn}\int_{2t}^{\infty}|X-Y|^{1-n}|F(Y)|dY\\
&\leq C\|F\|_{L^{p,\infty}(\rn\times (2t,\infty))}\bigg\|\big[|x-\cdot|^{2}+y_{n+1}^{2}\big]^{\frac{-(n-1)}{2}}\bigg\|_{L^{p',1}(\rn\times (2t,\infty))}.
\end{align}
Let $\varLambda(Y)=|Y|^{1-n}$, $Y=(y,y_{n+1})\in \rn\times (2t,\infty)$. Compute its distributional function 
\begin{align*}
d_{\varLambda}(s)=\bigg|\{(y,y_{n+1})\in \rn\times (2t,\infty)\}: |Y|<s^{\frac{1}{n-1}}\bigg|=c_nt^{n+1}\big((t^{n-1}s)^{-\frac{n+1}{n-1}}-1\big)
\end{align*}
so that by translation invariance, we find that $$\bigg\|\big[|x-y|^{2}+y_{n+1}^{2}\big]^{\frac{-(n-1)}{2}}\bigg\|_{L^{p',1}(\rn\times (2t,\infty))}= \displaystyle
\int^{(2t)^{1-n}}_{0}d^{\frac{1}{p'}}_{\varLambda}(s)ds=c_nt^{\frac{n+1}{p'}+1-n}$$
since $n+1>2p$. All together, this yields from \eqref{A4} the bound
\begin{align*}
|A_4|&\leq Ct^{2-\frac{n+1}{p}}\|F\|_{L^{p,\infty}(\rnp)}\leq Ct^{2-\frac{n+1}{p}}\|F\|_{\Y^{p}_{\infty}}.
\end{align*}
This shows \eqref{eq:weighted-sup}.
Next, we claim that
\begin{align}\label{eq:Cald-est}
\|\mathscr{G} F\big\|_{L^{\frac{(n+1)q}{n+1-q},\infty}(\rnp)}\leq C\|F\|_{\Y^{p}_{\infty}}.
\end{align}
Consider the even extension of a smooth function $F\in C^{\infty}_0(\rnp)$ to the whole space $\rnn$,
\begin{align*}
\widetilde{F}(X)=\begin{cases}
F(X)\hspace{0.2cm}\mbox{if}\hspace{0.1cm}X\in \rnp\\
F(x,-t)\hspace{0.2cm}\mbox{if}\hspace{0.1cm}X\in \mathbb{R}^{n+1}_{-}
\end{cases}
\end{align*}	
and compute 
\begin{align*}
I_{\alpha}\widetilde{F}(X)&=C(n,\alpha)\int_{\rnn}|X-Y|^{\alpha-(n+1)}\widetilde{F}(Y)dY\\
&=C(n,\alpha)\bigg[\int_{\rnn}\frac{\mathbf{1}_{\rnp}F(Y)dY}{|X-Y|^{(n+1)-\alpha}}+\int_{\rnn}\frac{\mathbf{1}_{\rnn_{-}}(Y)F(Y^{\ast})dY}{|X-Y|^{(n+1)-\alpha}}\bigg]\\
&=C(n,\alpha)\bigg[\int_{\rnp}\frac{F(Y)dY}{|X-Y|^{(n+1)-\alpha}}+\int_{\rnn_{-}}\frac{F(Y^{\ast})dY}{|X-Y|^{(n+1)-\alpha}}\bigg]\\
&=C(n,\alpha)\bigg[\int_{\rnp}\frac{F(Y)dY}{|X-Y|^{(n+1)-\alpha}}+\int_{\rnp}\frac{F(Y)dY}{|X-Y^{\ast}|^{(n+1)-\alpha}}\bigg]\\
&=\int_{\rnp}\bigg[\frac{C(n,\alpha)}{|X-Y|^{(n+1)-\alpha}}+\frac{C(n,\alpha)}{|X-Y^{\ast}|^{(n+1)-\alpha}}\bigg]F(Y)dY
\end{align*}
so that the Newtonian potential of $F$, $\mathscr{G}F(X)$ is the restriction to the half-space (up to a dimensional multiplicative constant) of the $\alpha$-Riesz potential defined on $\rnn$, that is,  $I_{\alpha}\widetilde{F}(X)\big|_{\rnp}=\mathscr{G} F(X)$. From this fact and by the mapping properties of the fractional operator between Lebesgue spaces, it follows that
\begin{align*}
\|\mathscr{G} F\|_{L^{d}(\rnp)}&
		\leq C\|I_{2} \widetilde{F}\|_{L^{d}(\rnn)}\\
		&\leq C\|\widetilde{F}\|_{L^{p}(\rnn)}\leq C\|F\|_{L^{p}(\rnp)}	
		\end{align*} 
provided we have $\dfrac{n+1}{d}=\dfrac{n+1}{q}-1$.
By a density argument and real interpolation, one gets the desired estimate
		\[\|\mathscr{G} F\|_{L^{\frac{(n+1)q}{n+1-q},\infty}(\rnp)}\leq C\|F\|_{L^{p,\infty}(\rnp)}\leq C\|F\|_{\Y^{p}_{\infty}}.\]
Regarding the gradient estimate 
\begin{align*}
\big\|\nabla \mathscr{G} F\big\|_{L^{\frac{(n+1)p}{n+1-p}}(\rnp)}
\leq  C\|F\|_{L^{p}(\rnp)},
\end{align*}
we exploit the pointwise decay property \eqref{eq:grad-bound} of the Green kernel. Indeed, we have 
\begin{equation*}
\nabla \mathscr{G}F(X)=\int_{\rnp}\nabla_{X}G(X,Y)F(Y)dY
\end{equation*}
in the sense of distributions and 
\begin{align}
\nonumber\big\|\nabla\mathscr{G}F\big\|_{L^{\frac{(n+1)p}{n+1-p}}(\rnp)}&=\bigg\|\int_{\rnp}\nabla_{X}G(X,Y)\varGamma(Y) dY\bigg\|_{L^{\frac{(n+1)p}{n+1-p}}(\rnp)}\\
\label{vertical-part-bd}&=\bigg\|\bigg(\int_{0}^{\infty}|\nabla GF(\cdot,t)|^{\frac{(n+1)p}{n+1-p}}dt\bigg)^{\frac{n+1-p}{(n+1)p}}\bigg\|_{L^{\frac{(n+1)p}{n+1-p}}(\rn)}.
\end{align}
Next, take $a>1$ and use Minkowski's and Young's inequalities to get  
\begin{align*}
\|\nabla\mathscr{G}F(x,\cdot)\|_{L^{\frac{(n+1)p}{n+1-p}}(0,\infty)}&=\bigg(\int_{0}^{\infty}|\nabla\mathscr{G}F(x,t)|^{\frac{(n+1)p}{n+1-p}}dt\bigg)^{\frac{n+1-p}{(n+1)p}}\\
&\leq C\bigg[\int_{0}^{\infty}\bigg(\int_{\rn}\int_{0}^{\infty}\dfrac{|F(Y)|}{[|x-y|^2+(t-y_{n+1})^2]^{\frac{n}{2}}}dY\bigg)^{\frac{(n+1)p}{n+1-p}}dt\bigg]^{\frac{n+1-p}{(n+1)p}}\\
&\leq C\bigg[\int_{0}^{\infty}\bigg(\int_{\rn}\mathcal{O}\ast |F|(x-y,t)dy \bigg)^{\frac{(n+1)p}{n+1-p}}dt\bigg]^{\frac{n+1-p}{(n+1)p}}\\
&\leq C\int_{\rn}|x-y|^{-n+\frac{1}{a}}\|F(y,\cdot)\|_{L^{p}(0,\infty)}dy
\end{align*} 
with $\dfrac{n+1-p}{(n+1)p}+1=\dfrac{1}{a}+\dfrac{1}{p}$ and $\mathcal{O}(x,t)=(|x|^2+t^2)^{-\frac{n}{2}}$ where we have used the fact that 
\begin{equation*}
\|\mathcal{O}(x,\cdot)\|_{L^{a}(0,\infty)}=|x|^{-n+\frac{1}{a}}\int_{0}^{\infty}(1+\sigma^2)^{\frac{-na}{2}}d\sigma,
\end{equation*}
the latter integral being finite since $an>1$. Therefore, the previous estimate implies the pointwise bound \begin{equation*}\|\nabla\mathscr{G}F(x,\cdot)\|_{L^{\frac{(n+1)p}{n+1-p}}(0,\infty)}\leq CI_{\frac{1}{a}}\|F(x,\cdot)\|_{L^{p}(0,\infty)}, \quad x\in \rn
\end{equation*}  so that the standard Riesz potential bounds between Lebesgue spaces yields in view of \eqref{vertical-part-bd}
\begin{align*}
\|\nabla\mathscr{G}F\|_{L^{\frac{(n+1)p}{n+1-p}}(\rnp)}\leq C\big\|I_{1/a}\|F(\cdot,t)\|_{L^{p}((0,\infty),dt)}\big\|_{L^{\frac{(n+1)p}{n+1-p}}(\rn)}\leq C\|F\|_{L^{p}(\rnp)} 
\end{align*}
since $\dfrac{n+1-p}{p(n+1)}=\dfrac{1}{p}-\dfrac{1/a}{n}$. The proof of \eqref{eq:grad-est-Green-kernel} follows the steps which have led to \eqref{eq:weighted-sup}, once again relying on the Green kernel bound \eqref{eq:grad-bound}.
This finishes the proof of Lemma \ref{lem:Green-pot-bound}.	
\end{proof}

Our next result deals with the mapping properties of the nonlinearity in Eq. \eqref{eq:main-eq}. 
\begin{lemma}\label{lem:mapping-prop}
Let $m>1$ and $1<q<n+1$. Put $\beta=\dfrac{(n+1)q}{(n+1-q)m}$ and consider the nonlinear map $N(u)=|u|^{m-1}u$. Then $N$ maps $\X^{q}_{\infty}$ onto $\Y^{\beta}_{\infty}$ $($resp. $\X^{q}$ onto $\Y^{\beta}$$)$ continuously. Moreover, we have the following bounds
 \begin{equation*}
 \|N(u)\|_{\Y^{\beta}_{\infty}}\leq C\|u\|^{m}_{\X^{q}_{\infty}}\hspace{0.2cm}\mbox{and}\hspace{0.2cm}
\|N(u)\|_{\Y^{\beta}}\leq C\|u\|^{m}_{\X^{q}}.
\end{equation*}
\end{lemma}
\begin{proof}
The proof follows directly from the definition of the spaces $\X^{q}_{\infty}$, $\Y^{\beta}_{\infty}$ and $\X^{q}$, $\Y^{\beta}$. The details are therefore omitted.	
\end{proof}
\section{Proofs of the main results}\label{S3}
In this section, we give detailed proofs of the results stated earlier in Section \ref{S-1} in the order of their appearance. 
\subsection{Proof of Theorem \ref{thm:1}} Let $q=\dfrac{(n+1)(m-1)}{m+1}$ and denote $q^{\ast}=\dfrac{(n+1)q}{n+1-q}$. Next, take $f\in L^{\frac{nq}{n+1},\infty}(\partial\rnp)$ and consider the operator $\mathscr{P}$ given by \[\mathscr{P}u=\mathscr{N}[bu|u|^{\eta-1}+f]+\mathscr{G}[a|u|^{m-1}u].\] In order to justify the statement $\textcolor{blue}{(i)}$ of Theorem \ref{thm:1}, it will suffice that the map $\mathscr{P}$ be contractive and self-mapping on a closed set in $\X^{q}_{\infty}$. Since the nonzero constant $a,b\in \mathbb{R}$ do not play a special role here, we simply omit them in this part. Let $u,v\in \X^{q}_{\infty}$ and write
\begin{align}\label{contraction}
\|\mathscr{P}u-\mathscr{P}v\|_{\X^{q}_{\infty}}&=\big\|\mathscr{N}[|u|^{\eta-1}u-|v|^{\eta-1}v]+\mathscr{G}[|u|^{m-1}u-|v|^{m-1}v]\big\|_{\X^{q}_{\infty}}\leq I+II
\end{align}  
where $I=\big\|\mathscr{N}[|u|^{\eta-1}u-|v|^{\eta-1}v]\big\|_{\X^{q}_{\infty}}$ and $II=\big\|\mathscr{G}[|u|^{m-1}u-|v|^{m-1}v]\big\|_{\X^{q}_{\infty}}$. Recall the Sobolev trace embedding
\begin{equation}\label{trace}
\mathscr{D}_{\infty}^{1,q}(\rnp)\hookrightarrow L^{\overline{q},\infty}(\partial\rnp),\quad \dfrac{n}{\overline{q}}=\dfrac{n+1}{q}-1, \quad 1<\overline{q}<q<\infty.    
\end{equation}
which follows by interpolating the inequality in \cite[Theorem 2]{Adam}; see also \cite[Corollary 1.4]{AL} for a more general result. One may invoke Proposition \ref{prop:linear-bounds} with $r_1=q^{\ast}$, $r_2=q$, $\mu=\infty$, $p=\dfrac{nq}{n+1}$ and exploit  the pointwise inequality
\begin{equation}\label{eq:point-ineq}\big||a_1|^{m-1}a_1-a_2|a_2|^{m-1}\big|\leq m|a_1-a_2|(|a_1|^{m-1}+|a_2|^{m-1})\quad
\end{equation}
 for all $a_1,a_2\in\mathbb{R}$ and $m>1$ to arrive at
\begin{align}\label{I}
\nonumber I&= \big\|\mathscr{N}[|u|^{\eta-1}u-|v|^{\eta-1}v]\big\|_{\X^{q}_{\infty}}\\
\nonumber&=\sup_{t>0}t^{\frac{n+1}{q}-1}\big\|\mathscr{N}[|u|^{\eta-1}u-|v|^{\eta-1}v]\big\|_{L^{\infty}(\rn)} +
\big\|\mathscr{N}[|u|^{\eta-1}u-|v|^{\eta-1}v]\big\|_{L^{q^{\ast},\infty}(\rnp)}+\\
\nonumber&\hspace{8.5cm}\big\|\nabla\mathscr{N}[|u|^{\eta-1}u-|v|^{\eta-1}v]\big\|_{L^{q,\infty}(\rnp)}\\
\nonumber&\leq C\big\|(u-v)(|u|^{\eta-1}+|v|^{\eta-1})\big\|_{L^{p,\infty}(\partial\rnp)}\\
\nonumber&\leq C\big\|u-v\big\|_{L^{\overline{q},\infty}(\partial\rnp)}\big\||u|^{\eta-1}+|v|^{\eta-1}\big\|_{L^{\frac{\overline{q}}{\eta-1},\infty}(\partial\rnp)}\\
\nonumber&\leq C\big\|u-v\big\|_{L^{\overline{q},\infty}(\partial\rnp)}(\|u\|^{\eta-1}_{L^{\overline{q},\infty}(\partial\rnp)}+\|v\|^{\eta-1}_{L^{\overline{q},\infty}(\partial\rnp)})\\
\nonumber&\leq C\|\nabla(u-v)\|_{L^{q,\infty}(\rnp)}\big(\|\nabla u\|^{\eta-1}_{L^{q,\infty}(\rnp)}+\|\nabla v\|^{\eta-1}_{L^{q,\infty}(\rnp)}\big)\\
I&\leq C\|u-v\|_{\X^{q}_{\infty}}(\|u\|^{\eta-1}_{\X^{q}_{\infty}}+\|v\|^{\eta-1}_{\X^{q}_{\infty}}).
\end{align}
Now, by utilizing Lemma \ref{lem:Green-pot-bound} \& Lemma \ref{lem:mapping-prop} together with H\"{o}lder inequality, we obtain
\begin{align}\label{II}
\nonumber II&= \big\|\mathscr{G}[u|u|^{m-1}-v|v|^{m-1}]\big\|_{\X^{q}_{\infty}}\leq C\big\|u|u|^{m-1}-v|v|^{m-1}\big\|_{\Y^{\frac{(n+1)q}{n+1+q}}_{\infty}}\\
\nonumber&\leq C\bigg[\sup_{t>0}t^{\frac{n+1}{q}+1}\big\|(u-v)[|u|^{m-1}+|v|^{m-1}](\cdot,t)\big\|_{L^{\infty}(\partial\rnp)}+\\
\nonumber&\hspace{7cm}\big\|(u-v)[|u|^{m-1}+|v|^{m-1}]\big\|_{L^{\frac{(n+1)q}{n+1+q},\infty}(\rnp)}\bigg]\\
\nonumber&\leq C\bigg[\sup_{t>0}t^{\frac{n+1}{q}-1}\|(u-v)(\cdot,t)\|_{L^{\infty}(\partial\rnp)}\sup_{t>0}t^{2}\big\|[|u|^{m-1}+|v|^{m-1}](\cdot,t)\big\|_{L^{\infty}(\partial\rnp)}+\\
\nonumber&\hspace{5,4cm}\|u-v\|_{L^{q^{\ast},\infty}(\rnp)}\big\|[|u|^{m-1}+|v|^{m-1}]\big\|_{L^{\frac{n+1}{2},\infty}(\rnp)}\bigg]\\
\nonumber&\leq C\bigg[\sup_{t>0}t^{\frac{n+1}{q}-1}\|(u-v)(\cdot,t)\|_{L^{\infty}(\partial\rnp)}\big((\sup_{t>0}t^{\frac{n+1}{q}}\|u(\cdot,t)\|_{L^{\infty}(\partial\rnp)})^{m-1}+\\
\nonumber&\hspace{.6cm}(\sup_{t>0}t^{\frac{n+1}{q}}\|v(\cdot,t)\|_{L^{\infty}(\partial\rnp)})^{m-1}\big)+\|u-v\|_{L^{q^{\ast},\infty}(\rnp)}\big(\|u\|^{m-1}_{L^{q^{\ast},\infty}(\rnp)}+\|v\|^{m-1}_{L^{q^{\ast},\infty}(\rnp)}\big)\bigg]\\
&\leq C\|u-v\|_{\X^{q}_{\infty}}(\|u\|^{m-1}_{\X^{q}_{\infty}}+\|v\|^{m-1}_{\X^{q}_{\infty}}).
\end{align}  
With $v=0$ in \eqref{contraction} taking into account \eqref{I} and \eqref{II}, one may apply Proposition \ref{prop:linear-bounds} once again to get
\begin{align*}
\|\mathscr{P}u\|_{\X^{q}_{\infty}}&\leq C(\|u\|^{m}_{\X^{q}_{\infty}}+\|u\|^{\eta}_{\X^{q}_{\infty}})+C'\|\mathscr{N}f\|_{\X^{q}_{\infty}}\\
&\leq C(\|u\|^{m}_{\X^{q}_{\infty}}+\|u\|^{\eta}_{\X^{q}_{\infty}})+C'\|f\|_{L^{\frac{nq}{n+1},\infty}(\rn)}
\end{align*}
for $f\in L^{\frac{nq}{n+1},\infty}(\rn)$. Take $\varepsilon>0$, set $\vartheta=2\varepsilon$ and let $B_{\vartheta}$ be the closed ball in $\X^{q}_{\infty}$ centered at the origin with radius $\vartheta$. If $\|f\|_{L^{\frac{nq}{n+1},\infty}(\partial\rnp)}\leq \varepsilon/C'$, then under the condition $C(\vartheta^{\eta}+\vartheta^{m})\leq \varepsilon$, a requirement which can be achieved provided $\varepsilon>0$ is taken sufficiently small, one deduces $\mathscr{P}:B_{\vartheta}\rightarrow B_{\vartheta}$ continuously and that $\mathscr{P}$ is a contraction from \eqref{I} and \eqref{II}. An application of the Banach fixed point Theorem infers the existence of a unique solution in a ball of $\X^q_{\infty}$.

Moving on, we prove the energy inequality \eqref{energy-ineq}. Assume that $f\in L^{\frac{2n}{n+2}}(\rn)$, $m=1+4/n$,  $n>2$. It is now clear (see Remark \ref{rmk:thm}) that if $\|f\|_{L^{\frac{nq}{n+1}}(\rn)}$ is sufficiently small, then there is a unique solution $u$ in $\X^{\frac{2(n+1)}{n+2}}$. In addition, it follows from \eqref{eq:est1} with $\kappa=1$ that \begin{equation*}
\displaystyle \sup_{t>0}t^{\frac{n}{2}+1}\|\nabla\mathscr{N}f(\cdot,t)\|_{L^{\infty}(\rn)}\leq C\|f\|_{L^{\frac{2n}{n+2}}}(\rn)    
\end{equation*} while 
\begin{equation*}
\displaystyle \sup_{t>0}t^{\frac{n}{2}+1}\|\nabla\mathscr{G}[|u|^{m-1}u](\cdot,t)\|_{L^{\infty}(\rn)}\leq C\|u^m\|_{L^{\frac{2(n+1)}{n+4}}(\rnp)}    
\end{equation*} is obtained by applying \eqref{eq:grad-est-Green-kernel} in Lemma \ref{lem:Green-pot-bound}. Thus, 
\begin{align*}
\int_{\rnp}t|\nabla u|^2dX&=\int_{\rnp}t|\nabla \mathscr{N}(f+b|u|^{\eta-1}u)|^2dX+\int_{\rnp}t|\nabla \mathscr{G}(a|u|^{m-1}u)|^2dX\\
&\leq C\bigg[(\sup_{t>0}t^{\frac{n}{2}+1}\|\nabla\mathscr{N}f(\cdot,t)\|_{L^{\infty}(\rn)})^{2/(n+2)}\int_{\rnp}|\nabla \mathscr{N}f|^{\frac{2(n+1)}{n+2}}dX+\\
&\hspace{0.5cm}(\sup_{t>0}t^{\frac{n}{2}+1}\|\nabla\mathscr{N}(b|u|^{\eta-1}u)(\cdot,t)\|_{L^{\infty}(\rn)})^{\frac{2}{n+2}}\int_{\rnp}|\nabla \mathscr{N}(|u|^{\eta-1}u)|^{\frac{2(n+1)}{n+2}}dX+\\
&\hspace{0.75cm}(\sup_{t>0}t^{\frac{n}{2}+1}\|\nabla\mathscr{G}[a|u|^{m-1}u](\cdot,t)\|_{L^{\infty}(\rn)})^{\frac{2}{n+2}}\int_{\rnp}|\nabla \mathscr{G}(|u|^{m-1}u)|^{\frac{2(n+1)}{n+2}}dX\bigg]\\
&\leq C\big(\|f\|^{2}_{L^{\frac{2n}{n+2}}(\rn)}+\|u\|^{2\eta}_{L^{2}(\rn)}+\|u\|^{2m}_{L^{\frac{2(n+1)}{n}}(\rnp)}\big)\\
&\leq C\big(\|f\|^{2}_{L^{\frac{2n}{n+2}}(\rn)}+\|\nabla u\|^{2\eta}_{L^{\frac{2(n+1)}{n+2}}(\rnp)}+\|u\|^{2m}_{L^{\frac{2(n+1)}{n}}(\rnp)}\big)\\
&\leq C\big(\|f\|^{2}_{L^{\frac{2n}{n+2}}(\rn)}+\|u\|^{2\eta}_{\X^{\frac{2(n+1)}{n+2}}}+\|u\|^{2m}_{\X^{\frac{2(n+1)}{n+2}}}\big)
\end{align*}
where to estimate the inequality before the last, we used $\mathscr{D}^{1,\frac{2(n+1)}{n+2}}(\rnp)\hookrightarrow L^{2}(\partial\rnp)$.
For $m=\dfrac{n+3}{n-1}$, one may repeat the previous argument using the boundedness properties of the operators $\mathscr{N}$ and $\mathscr{G}$ and exploiting the fact that $u$ lies in $\X^{q}$ to prove \eqref{energy-est-1}. This finishes the proof of Theorem \ref{thm:1}.    
\subsection{Proof of Theorem \ref{thm:2}}
As a direct consequence of the fixed point argument, the solution constructed in Theorem \ref{thm:1} may be realized as the limit in the space $\X^{q}_{\infty}$ of the sequence of Picard iterations given by
\[u_1(x,t)=\mathscr{N}f(x,t),\quad u_{j+1}=\mathscr{N}[u_j|u_j|^{\eta-1}]+\mathscr{G}[|u_j|^{m-1}u_j]+u_1,\hspace{0.1cm}j=1,2,...
\]  	
Once again we have dropped the coefficients here since their presence does not change the arguments. Let $q>1$ and $\varepsilon$ as before. Take $1<p_0<n$ and $0<\varepsilon_0<\varepsilon$ with $\varepsilon_0=\varepsilon_0(p_0)$. Given $f\in L^{p_0,\infty}(\rn)\cap L^{\frac{nq}{n+1},\infty}(\rn)$ satisfying $\|f\|_{L^{\frac{nq}{n+1},\infty}(\rn)}\leq \varepsilon_0$, we wish to show that  $(u_j)$ is a Cauchy sequence in $\X^{\frac{(n+1)p_0}{n}}_{\infty}$. Observe that this will produce the claim made in Theorem \ref{thm:3} since the limit of such a sequence solves  \eqref{eq:int-eq} and by uniqueness, it is nothing but the solution found in Theorem \ref{thm:1}. First observe that $u_j$ belongs to $\X^{\frac{(n+1)p_0}{n}}_{\infty}$ for each $j$. In fact, one has 
\begin{align}\label{eq:bound-Nf}\|u_1\|_{\X^{\frac{(n+1)p_0}{n}}_{\infty}}&=\sup_{t>0}t^{\frac{n}{p_0}-1}\|u_1(\cdot,t)\|_{L^{\infty}(\rn)}+\|u_1\|_{L^{\frac{(n+1)p_0}{n-p_0},\infty}(\rnp)}+\|\nabla u_1\|_{L^{\frac{(n+1)p_0}{n},\infty}(\rnp)}\\
\nonumber&\leq C_0\|f\|_{L^{p_0,\infty}(\rn)}
\end{align}
where we have applied Proposition \ref{prop:linear-bounds}. On the other hand, since $2(\eta-1)/(m-1)=1$, one may write $\frac{1}{p_0}=\frac{n-p_0}{np_0}+\frac{2(\eta-1)}{n(m-1)}$ and invoke Proposition \ref{prop:linear-bounds} once more,  H\"{o}lder's inequality and \eqref{trace} to arrive at
\begin{align}\label{eq:b1}
t^{-(\frac{n}{p_0}-1)}\nonumber\big|\mathscr{N}[u_1|u_1|^{\eta-1}]\big|&\leq C\big\|u_1|u_1|^{\eta-1}\big\|_{L^{p_0,\infty}(\rn)}\\
\nonumber&\leq C\|u_1\|_{L^{\frac{np_0}{n-p_0},\infty}(\rn)}\|u_1\|^{\eta-1}_{L^{\frac{n(m-1)}{2},\infty}(\rn)}\\
\nonumber&\leq C\|\nabla u_1\|_{L^{\frac{(n+1)p_0}{n},{\infty}}(\rnp)}\|\nabla u_1\|^{\eta-1}_{L^{q,{\infty}}(\rnp)}\\
 &\leq  C\|u_1\|_{\X^{\frac{(n+1)p_0}{n}}_{\infty}(\rnp)}\|u_1\|^{\eta-1}_{\X^{q}_{\infty}(\rnp)}\leq C\vartheta^{\eta-1}\|u_1\|_{\X^{\frac{(n+1)p_0}{n}}_{\infty}}.
\end{align}
With the same ingredients as before, one has
\begin{align}\label{b-2}
\nonumber\big\|\mathscr{N}[u_1|u_1|^{\eta-1}]\big\|_{L^{\frac{(n+1)p_0}{n-p_0},\infty}(\rnp)} 
&\leq C\big\|u_1|u_1|^{\eta-1}\big\|_{L^{p_0,\infty}(\rn)}\\
\nonumber&\leq C\|u_1\|_{L^{\frac{np_0}{n-p_0},\infty}(\rn)}\|u_1\|^{\eta-1}_{L^{\frac{n(m-1)}{2},\infty}(\rn)}\\
&\leq C\|u_1\|_{\X^{\frac{(n+1)p_0}{n}}_{\infty}}\|u_1\|_{\X^{q}_{\infty}}^{\eta-1}\leq C\vartheta^{\eta-1}\|u_1\|_{\X^{\frac{(n+1)p_0}{n}}_{\infty}}
\end{align}
and by using \eqref{eq:est2'} in Proposition \ref{prop:linear-bounds} with $r_2=\frac{(n+1)p_0}{n}$, $\mu=\infty$ one gets 
\begin{align}
\nonumber\big\|\nabla\mathscr{N}[u_1|u_1|^{\eta-1}]\big\|_{L^{\frac{(n+1)p_0}{n},\infty}(\rnp)}&\leq C\big\|u_1|u_1|^{\eta-1}\big\|_{L^{p_0,\infty}(\rn)}\leq C\vartheta^{\eta-1}\|u_1\|_{\X^{\frac{(n+1)p_0}{n}}_{\infty}}.
\end{align}
Next, by employing Lemmas \ref{lem:Green-pot-bound} and \ref{lem:mapping-prop}, we proceed similarly as above writing 
\begin{align}\label{b-3}
\nonumber\big\|\mathscr{G}[u_1|u_1|^{m-1}]\big\|_{\mathbf{X}^{\frac{(n+1)p_0}{n}}_{\infty}}&\leq C\big\|u_1|u_1|^{m-1}\big\|_{\Y^{(n+1)p_0/(n-p_0)}_{\infty}}\\
\nonumber&\leq C\big[\sup_{t>0}t^{\frac{n}{p_0}+1}\big\|u_1|u_1|^{m-1}(\cdot,t)\big\|_{L^{\infty}(\rn)}+\big\|u_1|u_1|^{m-1}\big\|_{L^{\frac{(n+1)p_0}{n+p_0},\infty}(\rnp)}\big]\\
\nonumber&\leq C\bigg[\sup_{t>0}t^{\frac{n}{p_0}-1}\|u_1\|_{L^{\infty}(\rn)}(\sup_{t>0}t^{\frac{n+1}{q}-1}\|u_1(\cdot,t)\|_{L^{\infty}(\partial\rnp)})^{m-1}+\\
\nonumber&\qquad\hspace{1.6cm} \|u_1\|_{L^{\frac{p_0(n+1)}{n-p_0},\infty}(\rnp)}\|u_1\|^{m-1}_{L^{q^{\ast},\infty}(\rnp)}
\bigg]\\
&\leq C\|u_1\|_{\X^{\frac{(n+1)p_0}{n}}_{\infty}}\|u\|^{m-1}_{\X^{q}_{\infty}}\leq C\vartheta^{m-1}\|u_1\|_{\X^{\frac{(n+1)p_0}{n}}_{\infty}}
\end{align}
where generalized H\"{o}lder's inequality was utilized to derive the third estimate.
From \eqref{eq:bound-Nf}, $u_1$ belongs to $\X^{\frac{(n+1)p_0}{n}}_{\infty}$. This fact, together with the estimates \eqref{eq:b1}, \eqref{b-2} and \eqref{b-3} yields the conclusion that $u_2$ also has a membership to $\X^{\frac{(n+1)p_0}{n}}_{\infty}$. Thus a simple induction argument permits us to deduce that $u_j\in \X^{\frac{(n+1)p_0}{n}}_{\infty}$ for each $j\in \mathbb{N}$. Assume that $C(\vartheta_0^{\eta-1}+\vartheta_0^{m-1})<1$,  $\vartheta_0=C_0\varepsilon_0$, $0<\varepsilon_0<\varepsilon$ where $C_0>0$ and $C>0$ are the constants appearing in \eqref{eq:bound-Nf} and in the sum \eqref{eq:b1}+\eqref{b-2}+\eqref{b-3}, respectively. Set $w_j=u_{j+1}-u_j$ and estimate the latter in $\X^{\frac{(n+1)p_0}{n}}_{\infty}$ using the pointwise inequality \eqref{eq:point-ineq} as follows
\begin{align*}
\|w_j\|_{\X^{\frac{(n+1)p_0}{n}}_{\infty}}&=\big\|\mathscr{N}[|u_j|^{\eta-1}u_j+|u_{j-1}|^{\eta-1}u_{j-1}]+\mathscr{G}[|u_j|^{m-1}u_{j}-|u_{j-1}|^{m-1}u_{j-1}]\big\|_{\X^{\frac{(n+1)p_0}{n}}_{\infty}}\\
&\leq C\big\|\mathscr{N}\big[w_{j-1}(|u_j|^{\eta-1}+|u_{j-1}|^{\eta-1})\big]\big\|_{\X^{\frac{(n+1)p_0}{n}}_{\infty}}+\\
&\hspace{6cm}C\big\|\mathscr{G}\big[w_{j-1}(|u_j|^{m-1}+|u_{j-1}|^{m-1})\big]\big\|_{\X^{\frac{(n+1)p_0}{n}}_{\infty}}\\
&\leq C\big[\|w_j\|_{\X^{\frac{(n+1)p_0}{n}}_{\infty}}(\|u_j\|^{\eta-1}_{\X^{q}_{\infty}}+
\|u_{j-1}\|^{\eta-1}_{\X^{q}_{\infty}})+\\
&\hspace{6.1cm}\|w_{j-1}\|_{\X^{\frac{(n+1)p_0}{n}}_{\infty}}(\|u_j\|^{m-1}_{\X^{q}_{\infty}}+
+\|u_{j-1}\|^{m-1}_{\X^{q}_{\infty}})\big]\\
&\leq C[2\vartheta_0^{\eta-1}+2\vartheta_0^{m-1}]\|w_{j-1}\|_{\X^{\frac{(n+1)p_0}{n}}_{\infty}}.
\end{align*}  
This inequality can be iterated finitely many times to obtain 
\begin{equation}\label{eq:iteration-est}
\|w_{j}\|_{\X^{\frac{(n+1)p_0}{n}}_{\infty}}\leq (2C)^{j-1}[\vartheta_0^{\eta-1}+\vartheta_0^{m-1}]^{j-1} \|w_{1}\|_{\X^{\frac{(n+1)p_0}{n}}_{\infty}}.
\end{equation}
Because $\|w_1\|_{\X^{\frac{(n+1)p_0}{n}}_{\infty}}<\infty$ and $2C[\vartheta_0^{\eta-1}+\vartheta_0^{m-1}]<1$, passing to the limit in \eqref{eq:iteration-est} as $j\rightarrow \infty$ gives $\|w_j\|_{\X^{\frac{(n+1)p_0}{n}}_{\infty}}\rightarrow 0$. This finishes the proof of Theorem \ref{thm:2}.

\subsection{Proof of Theorem \ref{thm:4}}

Let us reconsider Eq. \eqref{eq:int-eq} with the constants $a$ and $b$. Bear in mind that a solution of the latter is the limit in $\X^{q}_{\infty}$ of the sequence of approximations given by
\[u_1(x,t)=\mathscr{N}f(x,t);\hspace{0.1cm} u_{j+1}=\mathscr{N}[bu_j|u_j|^{\eta-1}]+\mathscr{G}[a|u_j|^{m-1}u_j]+u_1,\hspace{0.1cm}j=1,2,...
\]  
\textbf{Proof of Part \ref{pos-sol}}: Let $\Omega$ be a subset of $\rn$ with finite Lebesgue measure. Assume that $f>0$ in $\Omega$ in addition to $f$ being nonnegative in $\rn$. It is clear that $u_1$ is positive in $\rnp$ since the kernel $\mathcal{K}_t$ is positive for any $t>0$. Note that if $u_1\big|_{\partial\rnp}$ is nonnegative, then provided $b\geq0$, we have that $\mathscr{N}[bu_j|u_j|^{\eta-1}]\geq0$ in $\rnp$. Likewise $\mathscr{G}[au_1|u_1|^{m-1}]\geq0$ if $a\geq 0$ since the Green kernel $G(\cdot,\cdot)$ is positive. This shows that $u_2>0$ in $\rnp$. One may then proceed by induction to show that each element of the sequence $(u_j)$ is positive in $\rnp$. This sequence, however, converges in the space $\X^{q}_{\infty}$ (with $q>1$ as in Theorem \ref{thm:1}) to $u$ which satisfies 
\[u=\mathscr{N}[bu|u|^{\eta-1}]+\mathscr{G}[a|u|^{m-1}u]+u_1
\]    
in the sense of distributions.
Hence, by definition of $\X^{q}_{\infty}$, there is a subsequence of $(u_j)$ which converges pointwise almost everywhere to $u$. The conclusion immediately follows from the fact that almost everywhere convergence preserves positivity.\\
\textbf{Proof of Part \ref{rot-sym}}: Let $\mathscr{O}$ be the subfamily of rotations in $\rnp$ around the vertical axis. Given a point $(x,t)\in \rnp$, one may  write $\mathscr{O}(x,t)=(\mathscr{R}x,t)$ where $\mathscr{R}$ belongs to a group of rotations in $\rn$ preserving the origin. If $f$ is radial then $f(\mathscr{R}x)=f(x)$, $x\in \rn$. The kernels $\mathcal{K}_t$ and $G(x,t,\cdot)$ are both radial in the $x$-variable. This fact, together with the orthogonality of $\mathscr{R}$ justify our next calculations.
\begin{align*}
u_1[\mathscr{O}(x,t)]&=u_1(\mathscr{R}x,t)=\int_{\rn}\mathcal{K}_t(\mathscr{R}x-y)f(y)dy\\
&=\int_{\rn}\mathcal{K}_t(\mathscr{R}(x-\mathscr{R}^{-1}y))f(y)dy\\
&=\int_{\rn}\mathcal{K}_t(\mathscr{R}(x-z))f(\mathscr{R}z)dz\\
&=\int_{\rn}\mathcal{K}_t(x-z)f(z)dz=u_1(x,t)
\end{align*}
where we have made the change of variable $z=\mathscr{R}^{-1}y$. This implies that $u_1$ is radial. In the same spirit one may as well prove the identities \[\mathscr{N}[bu_1|u_1|^{\eta-1}](\mathscr{O}(x,t))=\mathscr{N}[bu_1|u_1|^{\eta-1}](x,t)\] and \[\mathscr{G}[au_1|u_1|^{m-1}](\mathscr{O}(x,t))=\mathscr{G}[au_1|u_1|^{m-1}](x,t)\] so that $u_2$ is radial in $\rnp$, and so is any $u_j$ through an induction procedure.  The convergence of $(u_j)$ to $u$ in $\X^{q}_{\infty}$ infers the convergence in $L^{\frac{(n+1)q}{n},\infty}(\rnp)$. Therefore, $u$ is rotationally symmetric around the vertical axis as convergence in the latter space preserves radial symmetry. The converse statement still flows from the radial symmetry properties of the kernels $\mathcal{K}_t$ and $G$. Indeed, if $u$ is rotationally symmetric, then so is $u_1=u-\mathscr{N}[bu|u|^{\eta-1}]-\mathscr{G}[au|u|^{m-1}]$. Its trace at the boundary, i.e. $I_1f$, the Riesz potential of order $1$ of $f$ is also radial. Necessarily, $f$ must be radial. Finally, we remark that the proof of Part \ref{rad-mono} is done exactly as above bearing in mind that "radially nonincreasing" is a property preserved by convolution.\\
\textbf{Proof of Part \ref{self-similar}:} Assume that $f\in L^{\frac{n(m-1)}{m+1},\infty}(\rn)$ and let $f_{\lambda}(x)=\lambda^{-\frac{m+1}{m-1}}f(x)$, $\lambda>0$ (e.g. $f(x)=|x|^{-\frac{m+1}{m-1}}$, $x\in \rn\setminus\{0\}$) and define $u_{\lambda}(x,t)=\lambda^{\frac{2}{m-1}}u(\lambda x,\lambda t).$ We compute 
\begin{align*}
[\mathscr{N}f]_{\lambda}&=\lambda^{\frac{2}{m-1}}\mathscr{N}f(\lambda x,\lambda t)
=\lambda^{\frac{2}{m-1}}\int_{\mathbb{R}^n}\mathcal{K}(\lambda x-y,\lambda t)f(y)dy\\
&=\lambda^{\frac{2}{m-1}}\lambda^{n}\int_{\mathbb{R}^n}\mathcal{K}(\lambda x-\lambda z,\lambda t)f(\lambda z)dz\\
&=\lambda^{\frac{2}{m-1}+1}\int_{\mathbb{R}^n}\mathcal{K}(x-z,t)f(\lambda z)dz\\
&=\int_{\mathbb{R}^n}\mathcal{K}(x-z,t)\big[\lambda^{\frac{m+1}{m-1}}f(\lambda z)\big]dz=\int_{\mathbb{R}^n}\mathcal{K}(x-z,t)f_{\lambda}(z)dz=\mathscr{N}(f_{\lambda}).
\end{align*} 
A similar computation will yield $[\mathscr{N}(b|u|^{\eta-1}u)]_{\lambda}=\mathscr{N}(b|u_{\lambda}|^{\eta}u_{\lambda})$
while the homogeneity property of the Green kernel, $G(X,Y)=\lambda^{n-1} G(\lambda X,\lambda Y)$ contributes to showing that \[[\mathscr{G}(a|u|^{m-1}u)]_{\lambda}=\mathscr{G}(a|u_{\lambda}|^{m}u_{\lambda}).\] Hence, $u_{\lambda}$ is a unique solution of Eq. \eqref{eq:int-eq}  
 whenever $f_{\lambda}$ has a small norm in $L^{\frac{n(m-1)}{m+1},\infty}(\rn)$ according to Theorem \ref{thm:1}. By hypothesis, $f$ and $f_{\lambda}$ have the same norm in $L^{\frac{n(m-1)}{m+1},\infty}(\rn)$. By uniqueness and in view of 
the relation $\|u_{\lambda}\|_{\X^{q}_{\infty}}=\|u\|_{\X^{q}_{\infty}}$, it follows that $u=u_{\lambda}$. The  proof of Theorem \ref{thm:4} is now complete.
\subsection*{Proof of Proposition \ref{prop:nonexistence-pos}} Let $a,b>0$ and $ m<(n+1)/(n-1)$. Assume that $f$ is positive and bounded below near the origin. By contradiction, suppose that Eq. \eqref{eq:int-eq} has a positive solution $u$ in $\X^{q}_{\infty}$. Then
\begin{align}\label{pointwise-nonexistence}
	\nonumber u(X)&=\mathscr{G}[au^{m}](X)+\int_{\rn}\dfrac{(bu^{\eta}+f)dy}{(|x-y|^2+t^2)^{\frac{n-1}{2}}}\geq c\int_{B_1(0)}\dfrac{f(y)dy}{(|x-y|^2+t^2)^{\frac{n-1}{2}}}\\
	&\geq C(|X|+1)^{-(n-1)}\quad \mbox{a.e.} \hspace{0.1cm}X\in \rnp.
	\end{align}
Let $R>0$, integrate the above inequality on the half-ball $B_R^{+}=B_R(0)\cap \rnp$, $B_R(0)\subset \rnn$ and multiply both sides of the resulting inequality by $|B_{R}^{+}|^{\frac{2}{(n+1)(m-1)}-1}$ to obtain
\begin{align*}
|B_{R}^{+}|^{\frac{2}{(n+1)(m-1)}-1}\bigg(\int_{B^+_R}\dfrac{dX}{(|X|+1)^{(n-1)}}\bigg)\leq C|B_{R}^{+}|^{\frac{2}{(n+1)(m-1)}-1}\int_{B^+_R}u(X)dX.
\end{align*}
For $R>0$ sufficiently large, this implies that \[R^{-(n-1)+\frac{2}{m-1}}\leq |B_{R}^{+}|^{\frac{2}{(n+1)(m-1)}-1}\|u\|_{L^{1}(B^+_R)}\leq C\|u\|_{L^{\frac{(n+1)q}{n+1-q},\infty}(\rnp)}.\]
Since $m<(n+1)/(n-1)$ by hypothesis, we find after passing to the limit in the previous inequality as $R\rightarrow \infty$ that $\|u\|_{L^{\frac{(n+1)q}{n+1-q},\infty}(\rnp)}$ blows up, thus contradicting the fact that $u$ belongs to $\X^q_{\infty}$. Next, suppose that $f$ is bounded below locally in $\rn$. Arguing as before, if $u$ is a positive solution of Problem \eqref{eq:int-eq} in $\X^{q}_{\infty}$, then  an analogue of the pointwise bound \eqref{pointwise-nonexistence} holds namely   
\begin{align*}
u(X)\geq \dfrac{CR^n}{(|X|+R)^{n-1}}\quad \mbox{a.e.}\hspace{0.2cm} X\in \rnp \hspace{0.1cm}\mbox{for any}\hspace{0.1cm}R>0.    
\end{align*} 
Integrate this inequality over $B_R^{+}$, $R$ chosen large enough to get 
\[R^{n+\frac{2}{m-1}-(n+1)}\int_{B_R^+}(|X|+R)^{-(n-1)}dX\leq C|B_R^+|^{\frac{2}{(n+1)(m-1)}-1}\int_{B^+_R}u(X)dX\leq C\|u\|_{L^{\frac{(n+1)q}{n+1-q},\infty}(\rnp)}\]
from which it follows (with $m=(n+1)/(n-1)$) that $R^n\leq C\|u\|_{\X^{q}_{\infty}}$ so that passing to the limit on both sides as $R\rightarrow \infty$ gives $C\|u\|_{\X^{q}_{\infty}}\geq \infty$ for some constant $C>0$ independent of $R$. This is a contradiction.  

\subsection*{Proof of Theorem \ref{thm:nonexistence-sol}} The proof is based on the test function method \cite{MP}. To this end we need a different notion of solutions. Let $f\in L^1_{loc}(\rn)$, call $u:\rnp\rightarrow \mathbb{R}$ a distributional solution to Eq. \eqref{eq:main-eq} (with positive nonlinearities) if it satisfies 
\begin{equation}
|u|^{\eta}\big|_{\partial\rnp}\in L^{1}_{loc}(\rn),\quad |u|^{m}\in L^{1}_{loc}(\rnp),\quad u\in L^{1}_{loc}(\rnp)    
\end{equation} 
and
\begin{equation}\label{eq:very-weak-sol}
\int_{\rnp}u\Delta \varphi dX-\int_{\rn}(b|u|^{\eta}+f)(x)\varphi (x,0)dx=a\int_{\rnp}|u|^{m}\varphi dX 
\end{equation} 
for all functions $\varphi\in C_0^{\infty}(\overline{\rnp})$ such that $\partial_t\varphi=0$ on $\partial\rnp$. We start with the following\\ 
\textbf{Claim:} A solution $u$ of Eq. \eqref{eq:int-eq} with $f$ locally integrable is a distributional solution of Eq. \eqref{eq:main-eq} in the sense of the above definition.\\
Let us momentarily defer the proof of this claim. Let $f$ be a locally integrable function and suppose that $\displaystyle\liminf_{|x|\rightarrow \infty}(1+|x|^{k})|f(x)|\geq L$ for some constant $L>0$ and some $0\leq k<\frac{m-1}{m+1}$. Next, let $\zeta \in C^{\infty}_0([0,\infty))$,  $0\leq \zeta \leq 1$ with
$$\zeta(\sigma)=\begin{cases}
1\hspace{0.2cm}\mbox{if}\hspace{0.2cm} \sigma\in [0,1]\\
0 \hspace{0.2cm}\mbox{if}\hspace{0.2cm} \sigma\in [2,\infty)
\end{cases}$$
and for $R>0$ large, introduce the rescaled functions 
$$\zeta^{1}_R(x)=\bigg[\zeta\bigg(\frac{|x|^{2}}{R}\bigg)\bigg]^{\frac{2m}{m-1}},\hspace{0.2cm} x\in \rn;\quad \zeta^{2}_R(t)=\bigg[\zeta\bigg(\frac{t^{2}}{R}\bigg)\bigg]^{\frac{2m}{m-1}},\hspace{0.2cm} t>0$$
and set $\zeta_R(x,t)=\zeta^{1}_R(x)\zeta^{2}_R(t)$. Clearly, $\zeta_R$ is smooth, has compact support in $\overline{\rnp}$ and $\partial_t\zeta_{R}=0$ on $\{t=0\}$. Arguing by contradiction, assume that Eq. \eqref{eq:int-eq} admits a solution $u$.  Then by the above claim (to be established later) $u$ is a solution of \eqref{eq:main-eq} in the sense of distributions and in particular, one has 
\begin{equation*}
\int_{\rn}u\Delta\zeta_R dX-\int_{\rn}(|u|^{\eta}+f)\zeta^{1}_R dx=\int_{\rnp}|u|^{m}\zeta_R dX\\
\end{equation*}
where we have assumed without any loss of generality that $a=b=1$. 
By applying Young's inequality, we obtain 
\begin{equation}\label{int1}
\int_{\rnp}|u|^m\zeta_R dX+\int_{\rn}(|u|^{\eta}+f)\zeta^{1}_R dx\leq \int_{\rnp}|u|^m\zeta_RdX+C_m\int_{\rnp}\zeta_R^{-\frac{1}{m-1}}|\Delta\zeta_R(X)|^{m'} dX
\end{equation}
which in turn implies that 
\begin{equation}\label{last-eq}
\displaystyle\int_{\rn}f\zeta^{1}_Rdx\leq C\int_{\rnp}\zeta_R^{-\frac{1}{m-1}}|\Delta\zeta_R(X)|^{m'} dX.
\end{equation}
A simple calculation shows that $\zeta_R$ satisfies the pointwise estimate 
\begin{align}\label{pointwise-est}
|\Delta\zeta_R (X)|\leq CR^{-1}\zeta_R^{1/m}(X)\big[\zeta^{1}_R(x)^{2}+\zeta^{2}_R(t)^{2}\big],\hspace{0.2cm}\hspace{0.2cm}X=(x,t)\in \rnp.    
\end{align}
In light of this, we have that 
\begin{align}\label{bound-test-fct}
\nonumber\mbox{R.H.S  of} \hspace{0.2cm}\eqref{last-eq}&\leq CR^{-m'}\bigg(\int^{\infty}_{0}\zeta_R^{2}(t)dt\bigg)|B_{(2R)^{1/2}}|+ CR^{-m'}R^{1/2}\int_{\rn}\zeta^{1}_R(x)dx\\
\nonumber&\leq CR^{-m'+\frac{n+1}{2}}\int_{0}^{\sqrt{2}}\zeta(\tau^{2})^{\frac{2m}{m-1}}d\tau+CR^{-m'+(n+1)/2}\int_{0}^{\sqrt{2}}\zeta(\tau^{2})^{\frac{2m}{m-1}}\tau^{n-1}d\tau\\
&\leq CR^{-m'+\frac{n+1}{2}}. 
\end{align}
By assumption, for any $\varepsilon>0$, there exits $R_0>0$ such that   
\[|f(x)|\geq (L-\varepsilon)(1+|x|^{k})^{-1}\hspace{0.2cm}\mbox{for all}\hspace{0.2cm} x\in \rn \hspace{0.2cm}\mbox{with} \hspace{0.2cm}|x|>R_0.\]
Hence, combining \eqref{last-eq} and \eqref{bound-test-fct} and for $R>R_0^{2}$ we have that
\begin{align*}
CR^{-m'+\frac{n+1}{2}}\geq \int_{\rn}f\zeta^{1}_Rdx\geq \int_{B_{R^{1/2}}}f\zeta^{1}_Rdx&=\int_{B_{R^{1/2}}}fdx\\
&=\int_{B_{R_0}}fdx+\int_{B_{R^{1/2}}\setminus B_{R_0}}fdx  \\
&\geq \int_{B_{R_0}}fdx+\int_{B_{R^{1/2}}}\dfrac{(L-\varepsilon)dx}{1+|x|^{k}}.
\end{align*}
We then deduce that
\begin{align*}
R^{-m'+\frac{n+1}{2}}&\geq C+C(1+R^{\frac{k}{2}})^{-1}R^{\frac{n}{2}}\geq  C(1+R^{\frac{n-k}{2}}).
\end{align*}
A fortiori, there must hold $\dfrac{n-k}{2}\leq -m'+\dfrac{n+1}{2}$ or $k\geq (m+1)/(m-1)$ which clearly yields a contradiction. To finish the proof of Theorem \ref{thm:nonexistence-sol}, one needs to justify the above claim.  
\begin{proof}[Proof of the Claim]
We proceed via integration by parts argument. Denote by 
$\langle\cdot,\cdot\rangle$ the duality pairing between $\mathbb{X}'$ and $\mathbb{X}=\{\varphi \in C_c^{\infty}(\overline{\rnp}):\partial_t \varphi =0 \hspace{0.2cm}\mbox{on}\hspace{0.2cm}\partial\rnp\}$. Let $u$ be a solution of Eq. \eqref{eq:int-eq} in $\X^q$, we clearly have $u\in L^1_{loc}(\rnp)$, $|u|^{\eta}\big|_{\partial\rnp}\in L^1_{loc}(\rn)$, $|u|^{m}\in L^1_{loc}(\rnp)$. For $\varphi\in \mathbb{X}$, we  have 
\begin{align*}
\langle u, \Delta \varphi\rangle=\nonumber\langle \mathscr{N}f,\Delta \varphi \rangle+\langle \mathscr{N}[b|u|^{\eta}],\Delta \varphi \rangle+\langle \mathscr{G}[a|u|^{m}],\Delta \varphi \rangle.
\end{align*}
Utilizing Fubini's theorem, we get that
\begin{equation}\label{int-part-Nf}
\langle \mathscr{N}f,\Delta \varphi \rangle=-\int_{\rnp}\mathscr{N}f(X)\Delta \varphi(X)dX=-\int_{\rn}\langle\nabla \mathcal{K}(\cdot-y),\nabla \varphi\rangle f(y)dy.
\end{equation} 
Next, we perform one integration by parts and use the harmonicity of $\mathcal{K}_t$ in $\rnp$ to obtain
\begin{align*}
\nonumber\langle\nabla\mathcal{K}(\cdot-y),\nabla\varphi\rangle&=\int_{\rnp}\nabla\mathcal{K}(X-y)\cdot\nabla\varphi(X)dX\\
\nonumber&=-\int_{\rnp}\Delta\mathcal{K}(X-y)\varphi(X)dX+\int_{\rn}\partial_{\nu}\mathcal{K}(x-y,0)\varphi(x,0)dx\\
\nonumber&=-\int_{\partial\rnp}\partial_{t}\mathcal{K}(x-y,0)\varphi(x,0)dx\\
\nonumber&=-\int_{\partial\rnp}\delta_y(x)\varphi(x,0)dx\\
&=-\varphi(y,0).
\end{align*} 
This, together with \eqref{int-part-Nf} gives $\langle \mathscr{N}f,\Delta \varphi \rangle=\displaystyle
\int_{\rn}f(y)\varphi(y,0) dy$. In the same vein, we also have $\langle \mathscr{N}[b|u|^{\eta}],\Delta \varphi \rangle=\displaystyle
\int_{\rn}b|u|^{\eta}(y)\varphi(y,0) dy$. 
Making use of Fubbini's theorem once more, we compute 
\begin{align}\label{int-part-Green-term}
\nonumber\big\langle \mathscr{G}[a|u|^{m}],\Delta \varphi \big\rangle&=\int_{\rnp}\mathscr{G}[a|u|^{m}](X)\Delta \varphi(X) dX\\
\nonumber&=\int_{\rnp}\int_{\rnp}G(X,Y)[a|u|^{m}](Y)dY\Delta\varphi(X)dX\\
&=-\int_{\rnp}\langle G(\cdot,Y),\Delta\varphi\rangle [a|u|^{m}](Y)dYdY.
\end{align}
But, 
\begin{align*}\label{green-int-parts}
\nonumber\langle \nabla G(\cdot,Y),\nabla \varphi\rangle&=\int_{\rnp}\nabla G(X,Y)\cdot\nabla\varphi(X)dX\\
\nonumber&=-\int_{\rnp}\Delta G(X,Y)\varphi(X)dX+\int_{\partial\rnp}\frac{\partial G}{\partial\nu}((x,0),Y)\varphi(x)dx\\
\nonumber&=-\int_{\rnp}\Delta G(X,Y)\varphi(X)dX\\
&=-\varphi(Y).
\end{align*}
so that $\big\langle \mathscr{G}[a|u|^{m}],\Delta \varphi \big\rangle=\displaystyle\int_{\rnp}a|u|^{m}\varphi dY$ from \eqref{int-part-Green-term}. Summarizing, $u$ fulfills the integral relation 
\begin{equation}
\int_{\rnp}u\Delta \varphi dX=\displaystyle\int_{\rnp}a|u|^{m}\varphi dX+\int_{\rn}(b|u|^{\eta}+f)\varphi(x,0) dx,
\end{equation} 
that is, $u$ is a distributional solution as defined in \eqref{eq:very-weak-sol}. This concludes the proof of Theorem \ref{thm:nonexistence-sol}.
\end{proof}

\subsection{Acknowledgements} The author wishes to thank his Ph.D advisor, Prof. Herbert Koch for his guidance and constant support.

\end{document}